%% file: main_VM.tex
\documentclass[11pt]{article}

\usepackage{graphicx}
\usepackage{latexsym,amsmath,amsfonts,amscd, amsthm, dsfont}
\usepackage{bm,color}
\usepackage{epsfig,verbatim,epstopdf,graphics}
\usepackage{subfigure}
\usepackage{changebar}
\usepackage{multirow}
\numberwithin{figure}{section}
\usepackage{xcolor}

	\usepackage[ruled,vlined]{algorithm2e}

	\usepackage{yhmath}
	\usepackage{booktabs} 
	\usepackage{tikz}
	\usepackage{verbatim}
	\usetikzlibrary{arrows,backgrounds,snakes,shapes}
	
	\numberwithin{equation}{section}

	\graphicspath{{./}{./figure/}}
	\allowdisplaybreaks
	
	\newtheoremstyle{plainNoItalics}{}{}{\normalfont}{}{\bfseries}{.}{ }{}
	
	\theoremstyle{plain}
	\newtheorem{thm}{Theorem}[section]
	
	\theoremstyle{plainNoItalics}

	\newtheorem{prop}[thm]{Proposition}
	\newtheorem{exa}[thm]{Example}

	\newcommand{\bx}{{\bf x}}

	\newcommand{\bv}{{\bf v}}
	\newcommand{\bw}{{\bf w}}
	\newcommand{\bB}{{\bf B}}
	\newcommand{\bE}{{\bf E}}
	\newcommand{\bJ}{{\bf J}}

	\newcommand{\bU}{{\bf U}}
	
	\newcommand{\bI}{{\bf I}}
	\newcommand{\bV}{{\bf V}}

	\newcommand{\bP}{{\bf P}}

	\newcommand{\br}{{\bf r}}

	\newcommand{\beq}{\begin{equation}}
		\newcommand{\eeq}{\end{equation}}
	\newcommand{\bit}{\begin{itemize}}
		\newcommand{\eit}{\end{itemize}}
	\newcommand{\be}{\begin{eqnarray}}
		\newcommand{\ee}{\end{eqnarray}}
	\newcommand{\beno}{\begin{eqnarray*}}
		\newcommand{\eeno}{\end{eqnarray*}}

	\newcommand{\Rmnum}[1]{\expandafter\@slowromancap\romannumeral #1@}
	\makeatother
	
	\usepackage{enumerate}

\begin{document}
\baselineskip=1.5pc

\vspace{.5in}

\input{title}

\input{intro}

\input{algorithm}

\input{AlgorithmNew}
\input{numerical}
\input{conclusion}

\newpage
\setcounter{page}{1} 
\bibliographystyle{siam}
\bibliography{refer,ref_guo,ref_cheng,ref_cheng_2}

\end{document}

%% file: title.tex
\begin{center}
{\bf
A Local Macroscopic Conservative (LoMaC) low rank tensor method for the Vlasov-Maxwell system
}
\end{center}

\vspace{.2in}
\centerline{
Shadi Heenatigala\footnote{
	Department of Mathematics and Statistics, Texas Tech University, Lubbock, TX, 70409. E-mail:                                                                                                                                       	shadiheenatigala92@gmail.com.
} 
		

and
 	
 Wei Guo\footnote{
	Department of Mathematics and Statistics, Texas Tech University, Lubbock, TX, 70409. E-mail:                                                                                                                                                                          
	weimath.guo@ttu.edu. Research is supported by NSF grant NSF-DMS-2111383 and Air Force Office of Scientific Research FA9550-22-1-0390.
} 
}

\bigskip
\noindent
{\bf Abstract.}

  The main computational challenges of solving the Vlasov-Maxwell (VM) system include the high dimensionality of the phase space, nonlinearity, inherent conservation properties, among others. In this paper, we develop a  novel Local Macroscopic Conservative (LoMaC) low rank tensor method for the VM system, as a continuation of our previous work (arXiv:2207.00518). The method takes advantage of the tensor friendly structure of the Vlasov equation and employs the low rank hierarchical Tucker decomposition to approximate the Vlasov solution in high dimensions. Hence, the curse of dimensionality can be mitigated. Furthermore, to realize the LoMaC property, the algorithm simultaneously evolves the  conservation laws of mass, momentum and energy  alongside the Vlasov equation using a  high order conservative method with the kinetic flux vector splitting. By a  conservative orthogonal projection, the low rank solution is guaranteed to have  the same macroscopic observables updated from the conservation laws. A collection of numerical tests on the VM system are presented to demonstrate the efficiency and efficacy of the proposed algorithm.

\vfill

\newpage

%% file: intro.tex
\section{Introduction}

In this paper, we are concerned with the Vlasov-Maxwell (VM) system, known as one of the most fundamental models in plasma physics. The VM system describes the dynamics of charged particles from the statistical mechanics viewpoint. In particular, the dimensionless single species non-relativistic VM system in 3D3V (three-dimensional physical space and three-dimensional velocity space) reads
\begin{align}
&	\partial_{t} f+\mathbf{v} \cdot \nabla_{\mathbf{x}} f+(\mathbf{E}+\mathbf{v} \times \mathbf{B}) \cdot \nabla_{\mathbf{v}} f=0\label{eq:vlasov}\\
&	\partial_{t} \mathbf{E}=\nabla_{\mathbf{x}} \times \mathbf{B}-\mathbf{J},\quad \partial_{t}\mathbf{B}=-\nabla_{\mathbf{x}} \times \mathbf{E} \label{eq:max1}\\
&	\nabla_{\mathbf{x}} \cdot \mathbf{E}=\rho-\rho_{i}, \quad \nabla_{\mathbf{x}} \cdot \mathbf{B}=0,\label{eq:max2}
\end{align}
where the unknown function $f(\bx,\bv,t)$ of the Vlasov equation \eqref{eq:vlasov} is the probability distribution function of electrons on the phase space $\Omega_{\bx}\times \Omega_{\bv}$, $\bx=(x_1,x_2,x_3)$ and $\bv=(v_1,v_2,v_3)$.  The electromagnetic fields $\bE=(E_1,E_2,E_3)$ and $\bB=(B_1,B_2,B_3)$ are determined from Maxwell's equations \eqref{eq:max1}-\eqref{eq:max2}.  The  density $\rho$ and current density $\bJ$ are given by 
$$\rho(\bx,t)=\int_{\Omega_\bv} f(\bx,\bv,t)d\bv,\quad \bJ(\bx,t)=\int_{\Omega_\bv} \bv f(\bx,\bv,t)d\bv.$$
The ion is assumed to be fixed as a uniform neutralizing background with density $\rho_i$. Furthermore, by taking the first few moments of the VM system, a set of macroscopic equations can be derived
\begin{align}
	\partial_{t} \rho + \nabla_\bx \cdot \bJ &= 0\label{eq:mass}\\
	\partial_{t} \bP +\nabla_{\bx} \cdot \bm{\sigma}&= \nabla_{\bx} \cdot
	\left(\bE\otimes\bE + \bB\otimes\bB -\frac12(|\bE|^2+|\bB|^2)\bI \right) + \bE\rho_i \label{eq:mom}\\
	\partial_{t} e +\nabla_{\bx} \cdot \mathbf{Q}& = \nabla_{\bx} \cdot(\bE\times\bB).\label{eq:ener} 
\end{align}
Here $\bP(\bx,t) = \bJ(\bx,t) + \bE\times\bB$ denotes the momentum density. $e(\bx,t)=\kappa(\bx,t) + \frac12|\bE|^2 + \frac12|\bB|^2$ denotes the energy density as the sum of  the kinetic energy density $\kappa(\bx,t)=\int_{\Omega_\bv} \frac12|\bv|^2 f(\bx,\bv,t)d\bv$, electric energy density $\frac12 |\bE|^2$, and the magnetic energy density $\frac12|\bB|^2$.  $ \bm{\sigma}(\bx,t)=\int_{\Omega_{\bv}}\bv\otimes\bv f(\bx, \bv,t) d \bv$ and $\mathbf{Q}(\bx,t) =\frac12\int_{\Omega_{\bv}}\bv|\bv|^2 f(\bx, \bv,t) d \bv$  are the flux terms, and $\bI$ is the identity matrix. Equations \eqref{eq:mass}-\eqref{eq:ener} correspond to the conservation laws of mass, momentum and energy, respectively. Note that the total momentum of the single species VM system is conserved only if $\int_{\Omega_\bx} \bE\, d\bx = 0$, which can be seen from \eqref{eq:mom}. The numerical challenges include the high dimensionality of phase space,  inherent conservation properties with respect to the macroscopic equations \eqref{eq:mass}-\eqref{eq:ener}, nonlinearity, multi-scale features, among others.
In this paper, we propose a novel deterministic  VM solver that is resistant to the curse of dimensionality and meanwhile can conserve the physical invariants locally at the discrete level.

The deterministic kinetic schemes are designed in the grid-based Eulerian or semi-Lagrangian (SL) fashion, and they are free of statistical noise and can attain high order accuracy conveniently, as opposed to particle methods, see e.g., \cite{filbet2003comparison,sircombe2009valis,cheng2014discontinuous,crouseilles2015hamiltonian,li2019solving}. The development of deterministic Vlasov solvers has elicited lots of attention, since they are capable of correctly capturing many  solution structures of practical interest, such as small-scale structures in low density regions or at the tail of the distribution, while particle methods are often not efficient.  Meanwhile, a standard deterministic discretization of the Vlasov equation will suffer the curse of dimensionality, making the method prohibitively expensive. To overcome the difficulty, some dimension reduction techniques are proposed for the deterministic approach  to mitigate the  curse of dimensionality, such as the sparse grid approach \cite{kormann2016sparse,guo2016sparse1}.  In this work, we are concerned with another effective tool for dimension reduction, namely low rank tensor decomposition. Such an approach has been developed with great success for approximating functions in high dimensions \cite{hackbusch2012tensor}, as well as for simulating high-dimensional partial differential equations (PDEs) including the Vlasov system using various types of tensor formats, see e.g., \cite{kormann2015semi,ehrlacher2017dynamical}.  In addition, a class of dynamical low rank methods are  proposed in \cite{einkemmer2018low,einkemmer2020low,dektor2020dynamically,einkemmer2021mass,cassini2022efficient} for which the dynamical low rank approximation of the Vlasov solution is evolved on the low rank manifold with a tangent space projection.

  Recently, we initialized a line of research for the development of low rank tensor methods for the deterministic Vlasov-Poisson (VP) simulations. In particular, with a low rank tensor representation of the solution, the proposed method takes advantage of the tensor friendly structure of the Vlasov equation and adds new basis by applying the well-established high order finite difference upwind method coupled with the strong-stability-preserving (SSP) multi-step time discretizations \cite{gottlieb2011strong}, followed by an SVD-type truncation to remove redundancy. Thanks to the hierarchical Tucker (HT) decomposition \cite{hackbusch2009new,grasedyck2010hierarchical}, the low rank algorithm attains a storage complexity that is linearly scaled with the dimension. On the other hand, there is a drawback of the proposed method that the SVD-type truncation procedure would destroy any conservation properties of the discretization and may result in unphysical behaviors of the numerical solution, especially in long term simulations. To overcome the difficulty, inspired by the work in \cite{einkemmer2021mass}, we developed a novel conservative truncation algorithm  that can remove the redundancy in the solution for computational efficiency and in the meantime, preserve the mass, momentum and kinetic energy densities of the pre-compressed solution \cite{guo2022conservative}. More recently, to achieve the desired local energy conservation at the discrete level, we developed  a class of Local Macroscopic Conservative (LoMaC) low rank tensor methods. The key idea is that we simultaneously update macroscopic conservation laws (e.g., \eqref{eq:mass}-\eqref{eq:ener}) alongside the VP system. In particular, the low rank Vlasov solution is used to construct numerical fluxes to update the macroscopic densities via the kinetic flux vector splitting (KFVS) \cite{mandal1994kinetic, xu1995gas}. Meanwhile, with the conservative truncation algorithm, after the adding basis step the pre-compressed  Vlasov solution is orthogonally projected onto the reference subspace defined by macroscopic densities updated from the conservation laws, followed by a conservative weighted truncation step to remove the redundancy. Hence, the macroscopic fluid solver and the low rank kinetic solver are implemented alongside with each other in a self-consistent fashion, and to the best of our knowledge, this is the first explicit low rank VP solver that achieves local energy conservation at the fully discrete level with theoretical guarantee. Note that the LoMaC property automatically leads to  global conservation of the total mass, momentum and energy of the system and is of much more practical interest to plasma simulations \cite{kraus2017gempic}.

 In this paper, we extend the proposed LoMaC VP solver to the VM system \eqref{eq:vlasov}-\eqref{eq:max2} by simultaneously considering the associated with macroscopic conservation laws \eqref{eq:mass}-\eqref{eq:ener}. In particular, the Vlasov solution is represented in the HT tensor format to mitigate the curse of dimensionality, and the spectral collocation method \cite{gottlieb2001spectral,hesthaven2007spectral} is employed for discretizing differential operators in the adding basis procedure. Furthermore, the proposed method achieves the LoMaC property by simultaneously  working with the low rank VM solver and the macroscopic fluid solver with explicit time marching. Noteworthy,  as with the LoMaC VP solver \cite{guo2022local}, the energy is conserved locally at the fully discrete level in the explicit manner, as opposed to existing energy conserving VM solvers which require implicit symplectic time integrators \cite{cheng2014energy1,cheng2014energy}.

This paper is organized as follows. In Section 2, we first review the low rank tensor approach for the 1D2V VM system in Section 2.1, then we introduce the key ingredients including the conservative low rank truncation algorithm in Section 2.2, together with the LoMaC algorithm with simultaneous update of macroscopic conservation laws using KFVS in Section 2.3. In Section 3, we present an extensive set of benchmark tests to demonstrate the effectiveness and the conservation properties of the proposed low rank tensor algorithm.  We conclude the main contributions of the paper and comment on future research directions in Section 4.

%% file: algorithm.tex

\section{A low rank tensor approach for the Vlasov-Maxwell system with local conservation}

We consider the following reduced 1D2V VM system 
\beq
\frac{\partial f}{\partial t}+v_{1} \frac{\partial f}{\partial x_{1}}+\left(E_{1}+v_{2} B_{3}\right) \frac{\partial f}{\partial v_{1}}+\left(E_{2}-v_{1} B_{3}\right) \frac{\partial f}{\partial v_{2}}=0,
\label{vlasov1}
\eeq
\begin{equation}
\label{eq:max11d}
		\frac{\partial E_{1}}{\partial t}=-J_1,\quad
		\frac{\partial E_{2}}{\partial t}=-\frac{\partial B_{3}}{\partial x_{1}}-J_2,\quad
		 \frac{\partial B_{3}}{\partial t}=-\frac{\partial E_{2}}{\partial x_{1}},
\end{equation}
\begin{equation}
\frac{\partial E_{1}}{\partial x_{1}}=\rho-\rho_i, \label{eq:gauss1d}
\end{equation}
where the unknown  $f(x_1,v_1,v_2,t)$ denotes the probability distribution function of electrons in the reduced phase space. 
The associated macroscopic conservation laws become 
\begin{align} 
	&\frac{\partial \rho}{\partial t} + \frac{\partial J_1}{\partial x_1} = 0,\label{eq:mass1d}\\
	& \frac{\partial P_1}{\partial t} + \frac{\partial \sigma_1}{\partial x_1} =\frac12\frac{\partial}{\partial x_1}\left(|E_1|^2 - |E_2|^2 - |B_3|^2\right) + \rho_iE_1, \label{eq:mom11d}\\
	& \frac{\partial P_2}{\partial t} + \frac{\partial \sigma_2}{\partial x_1} = \frac{\partial\left(E_1E_2\right)}{\partial x_1} + \rho_iE_2, \label{eq:mom21d} \\
	& \frac{\partial e}{\partial t} + \frac{\partial Q_1}{\partial x_1} =\frac{\partial\left(E_2B_3\right)}{\partial x_1}, \label{eq:ener1d}
\end{align}
where $P_1 = J_1 + E_2B_3$, $P_2 = J_2 - E_1B_3$, $\sigma_1 = \int v_1^2f dv_1dv_2$, $\sigma_2 = \int v_1v_2f dv_1dv_2$, and $Q_1 =\frac12\int v_1(v_1^2 +v_2^2 ) f dv_1dv_2 $.

\subsection{A low rank tensor approach for the Vlasov-Maxwell system}

We assume a spectral collocation discretization of $f(x_1,v_1,v_2)$, denoted by $\mathbf{f}$, on a truncated 1D2V domain of $[x_{1,\min}, x_{1,\max}] \times [-v_{1,\max}, v_{1,\max}]\times [-v_{2,\max}, v_{2,\max}]$ with uniform tensor product $N_{x_1} \times N_{v_1}\times N_{v_2}$ grid points 
\beq
\label{eq: x_grid}
x_{1,\text{grid}}: \quad x_{1,\min}=x_{1,1}< \cdots < x_{1,i} < \cdots < x_{1,N_{x_1}} = x_{1,\max}, 
\eeq
\beq
\label{eq: v_grid}
v_{1,\text{grid}}: \quad -v_{1,\max}=v_{1,1}< \cdots < v_{1,j} <\cdots < v_{1,N_{v_1}} = v_{1,\max},
\eeq
\beq
\label{eq: v_grid2}
v_{2,\text{grid}}: \quad -v_{2,\max}=v_{2,1}< \cdots < v_{2,j} <\cdots < v_{2,N_{v_2}} = v_{2,\max}.
\eeq
Denote $h_x$, $h_{v_1}$ $h_{v_2}$ as the mesh sizes in $x$, $v_1$ and $v_2$ directions, respectively, and denote $N=\max\{N_{x_1},N_{v_1},N_{v_2}\}$.

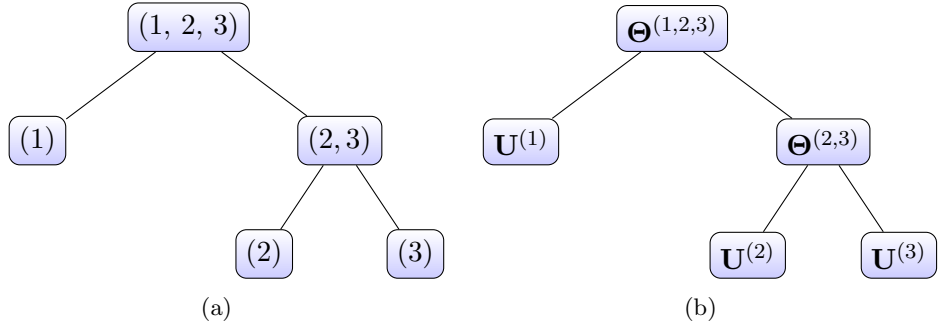
\begin{figure}
\centering
\subfigure[]{
 \begin{tikzpicture}[
     level/.style={sibling distance=40mm/#1},
  every node/.style = {shape=rectangle, rounded corners,
    draw, align=center,
    top color=white, bottom color=blue!20}
    ]
  \node {$(1,\,2,\,3)$}
    child { node {$(1)$} 
    }
    child { node {$(2,3)$}
    	child{ node{$(2)$}}
	child{ node{$(3)$}} 
	};
\end{tikzpicture}}\quad
\subfigure[]{
\begin{tikzpicture}[
  every node/.style = {shape=rectangle, rounded corners,
    draw, 
    top color=white, bottom color=blue!20},
    level/.style={sibling distance=40mm/#1}
    ]
  \node {$\bm{\Theta}^{(1,2,3)}$}
      child { node {$\bU^{(1)}$} 
    }
    child { node {$\bm{\Theta}^{(2,3)}$}
    	child{ node{$\bU^{(2)}$}}
	child{ node{$\bU^{(3)}$}} 
	};
\end{tikzpicture}}
 \caption{Dimension tree $\mathcal{T}$ and associated data layout to express third order tensors in the HT format. }
\label{fig:dimtree1}
\end{figure}

Directly working with the tensor product full grid representation of $\bf f$ suffers the curse of dimensionality, and  the proposed low rank tensor approach  is designed based on the low rank tensor decomposition of $\mathbf{f}$ in the HT format to mitigate the curse of dimensionality. The HT format is characterized by a dimension tree $\mathcal{T}$, frames at leaf nodes and transfer tensors at non-leaf nodes. In particular,  we denote by the dimension index $D:=(1,2,3)$    and employ the dimension tree $\mathcal{T}$ depicted in Figure \ref{fig:dimtree1} (a).  $\mathcal{T}$ is a binary tree and has $D$ as the root node and has $(1),( 2), (3)$ as the leaf nodes corresponding to dimension $x_1$, $v_1$, and $v_2$, respectively. In Figure \ref{fig:dimtree1} (b), the associated data structure is given: at a non-leaf node, a third order transfer tensor $\bm{\Theta}$ is stored, and at a leaf node a frame $\bU$ is stored instead. For example, ${\bf f}$ can be expressed as
 \beq
\label{eq:htd_f_nested_0}
{\bf f} = \sum_{l_{1}=1}^{r_{1}}{\sum_{l_{23}=1}^{r_{23}}} \bm{\Theta}^{(1,2, 3)}_{l_{1},l_{23}}\bU_{l_{1}}^{(1)}\otimes \bU_{l_{23}}^{(2, 3)},
\eeq
  with 
\beq
\label{eq: U34_htd}
\bU_{l_{23}}^{(2, 3)} = \sum_{l_{2}=1}^{r_{2}}{\sum_{l_{3}=1}^{r_{3}}} \bm{\Theta}^{(2, 3)}_{l_{2},l_{3}, l_{23}}\bU_{l_{2}}^{(2)}\otimes \bU_{l_{3}}^{(3)},\quad l_{23} =1,\ldots,r_{23}. 
\eeq
Denote $ \br = \{r_\alpha\}_{\alpha\in\mathcal{T}}$ as the hierarchical ranks. The storage of the HT format scales as $\mathcal{O}(r^3+ r(N_{x_1}+N_{v_1}+N_{v_2}))$, where $r=\max \br$. If $r$ is reasonably low, then the HT format avoids the curse of dimensionality. 
 
 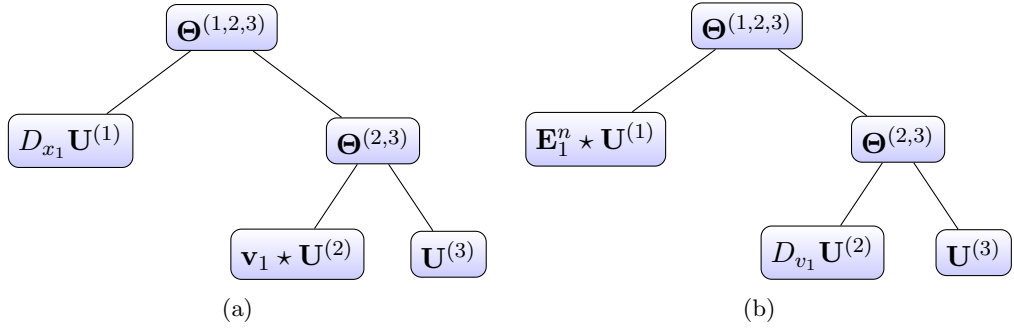
\begin{figure}
\centering
\subfigure[]{
\begin{tikzpicture}[
  every node/.style = {shape=rectangle, rounded corners,
    draw, 
    top color=white, bottom color=blue!20},
    level/.style={sibling distance=40mm/#1}
    ]
  \node {$\bm{\Theta}^{(1,2,3)}$}
      child { node {$D_{x_1}\bU^{(1)}$} 
    }
    child { node {$\bm{\Theta}^{(2,3)}$}
    	child{ node{$\bv_1\star\bU^{(2)}$}}
	child{ node{$\bU^{(3)}$}} 
	};
\end{tikzpicture}}\quad
\subfigure[]{
\begin{tikzpicture}[
  every node/.style = {shape=rectangle, rounded corners,
    draw, 
    top color=white, bottom color=blue!20},
    level/.style={sibling distance=40mm/#1}
    ]
  \node {$\bm{\Theta}^{(1,2,3)}$}
      child { node {$\bE_1^n\star\bU^{(1)}$} 
    }
    child { node {$\bm{\Theta}^{(2,3)}$}
    	child{ node{$D_{v_1}\bU^{(2)}$}}
	child{ node{$\bU^{(3)}$}} 
	};
\end{tikzpicture}}
 \caption{Data layout to express $v_1\partial_{x_1}(f^n)$ and $E_1^n\partial_{v_1}(f^n)$  in the HT format. }
\label{fig:dimtree2}
\end{figure}


The proposed low rank  tensor approach adaptively updates low rank basis and transfer tensors by two steps: an adding basis step by the spectral collocation method and a removing basis step by hierarchical high order SVD (HOSVD) truncation \cite{grasedyck2010hierarchical}. We apply a second order SSP multi-step temporal discretization  to illustrate the main idea.   $\mathbf{f}^n$ denotes the low rank solution tensor in the HT format at the time level $t^n$.
\begin{enumerate}
\item {\em Add basis and obtain an intermediate solution ${\bf f}^{n+1, *}$.} 
A second order multi-step discretization of time derivative in \eqref{vlasov1} gives
\begin{equation}
\label{eq: fn3}
{f}^{n+1, *} = \frac14{f}^{n-2}+\frac34 {f}^{n}- \frac32\Delta t \left(v_1 \partial_{x_1} ({f}^n) + E^n_1 \partial_{v_1} ({f}^n) + v_2B^n_3 \partial_{v_1} ({f}^n) + E^n_2 \partial_{v_2} ({f}^n) -v_1B^n_3 \partial_{v_2} ({f}^n) \right).
\end{equation}
Here the electric and magnetic fields $E^n_1$, $E^n_2$, $B^n_3$ are solved from Maxwell's equations. Thanks to the tensor friendly form of the Vlasov equation, each term in $f^{n+1,*}$ can be approximated by  HT tensors. For example, $v_1 \partial_{x_1} ({f}^n)$ and $E^n_1 \partial_{v_1} ({f}^n)$ are expressed in the HT format with the same dimension tree $\mathcal{T}$ of $\mathbf{f}^{n}$ together with the data layout given in Figure \ref{fig:dimtree2}, where, with a slight abuse of notation, $\bv_1 \in\mathbb{R}^{N_{v_1}}$ denotes the coordinates of $v_{1,\text{grid}}$ introduced in \eqref{eq: v_grid}. The discrete differentiation operators $D_{x_1}$ and $D_{v_1}$ from the spectral collocation method implemented with FFT  are applied to the basis in the dimension-by-dimension fashion, and $\star$ denotes an element-wise multiplication operation. 

\item {\em Remove basis of ${\bf f}^{n+1, *}$ and update ${\bf f}^{n+1}$.} Since the sizes of bases  and transfer tensors have increased in the adding basis step,  we perform a hierarchical HOSVD truncation to remove redundancy with prescribed threshold $\varepsilon$. Note that any conservation properties of mass, momentum or energy are lost after truncation. The removing basis step costs $\mathcal{O}(r^2 N + r^4)$ where $N=\max(N_{x_1}, N_{v_1},N_{v_2})$, see \cite{grasedyck2010hierarchical}.   
\end{enumerate}

\subsection{Conservative truncation} 
 In \cite{guo2022conservative}, we proposed a conservative low rank truncation algorithm
for preservation of mass, momentum and kinetic energy densities inspired by the work in \cite{einkemmer2021mass}. The key idea is to first project the pre-compressed solution ${\bf f}$ to a subspace 
 \beq
 \label{eq: 2d2v_subspa}
 \mathcal{N}
 = \text{span}\{{\bf 1}_{v_1} \otimes {\bf 1}_{v_2}, {\bf v}_1\otimes {\bf 1}_{v_2},  {\bf 1}_{v_1}\otimes{\bf v}_2, {\bf v}_1^2\otimes {\bf 1}_{v_2}+{\bf 1}_{v_1}\otimes{\bf v}_2^2\},
 \eeq
where  ${\bf 1}_{v_k}\in \mathbb{R}^{N_{v_k}}$ is the vector of all ones, $\bv_k \in\mathbb{R}^{N_{v_k}}$ denotes the coordinates of $v_{k,\text{grid}}$, and $\bv_k^2$ $\in \mathbb{R}^{N_{v_k}}$ is the element-wise square of $\bv_{k}$, $k=1,\, 2$. To ensure proper decay of the projected function as $v_1,v_2 \to \infty$, 
a weight function $w(v_1,v_2)$ is introduced with exponential decay.  One such example is $w(v_1,v_2) :=w^{(1)}(v_1)w^{(2)}(v_2) = \exp(-\frac{v_1^2}{2}) \exp(-\frac{v_2^2}{2})$. With  the weight function, a scaling and re-scaling procedure is needed for the projection step, as well as for the truncation step. 

The following definitions are introduced to formulate the algorithm:
\bit
\item Standard $l^2$ inner product: 
\beq
\label{eq: vm_inner_prod_2_d}
\langle {\bf f},  {\bf g} \rangle = h_{v_1}h_{v_2}\sum_{j_1=1}^{N_{v_1}}\sum_{j_2=1}^{N_{v_2}}  f_{j_1,j_2} g_{j_1,j_2}.
\eeq 

\item Weighted inner products and the associated norms:
\beq
 \label{eq: 1d_inner_w}
  \langle {\bf f},  {\bf g}\rangle_{{\bf w}^{(k)}} = h_{v_k} \sum_{j_m=1}^{N_{v_k}}   f_{j_m} g_{j_k} w^{(k)}_{j_k}, \quad \|\mathbf{f}\|_{\bw^{(k)}} = \sqrt{ \langle f, f\rangle_{{\bf w}^{(k)}}},
 \eeq
\beq
 \label{eq: vm_inner_w}
  \langle {\bf f},  {\bf g}\rangle_{{\bf w}} = h_{v_1}h_{v_2} \sum_{j_1=1}^{N_{v_1}} \sum_{j_2=1}^{N_{v_2}}  f_{j_1,j_2} g_{j_1,j_2} w^{(1)}_{j_1} w^{(2)}_{j_2}, \quad \|\mathbf{f}\|_\bw = \sqrt{ \langle f, f\rangle_{{\bf w}}},
 \eeq
where ${\bf w}^{(k)}\in \mathbb{R}^{N_{v_k}}$ with $w^{(k)}_j = w^{(k)}(v_j)$, $k=1,\, 2$, and $\bw = {\bf w}^{(1)}\otimes {\bf w}^{(2)}\in  \mathbb{R}^{N_{v_1}\times N_{v_2}}$. The definition is in analog to the weighted inner products at the continuous level, and we let
$
l^2_{\bf w} = \{{\bf f}\in\mathbb{R}^{N_{v_1}\times N_{v_2}}: \|{\bf f}\|_{\bf w} < \infty\}.
$
\eit 
Consider the subspace $\mathcal{N}\subset l^2_\bw$, a conservative low rank truncation of a numerical solution ${\bf f}$ written in the low rank HT tensor format can be obtained from the procedure below. 

\begin{enumerate}
\item {\bf Compute the discrete macroscopic densities.}
The discrete macroscopic charge, current and kinetic energy densities of ${\bf f}$ are given by
\begin{align}
\left(\begin{array}{l}
{\boldsymbol\rho}\\
{\bf J}_1\\
{\bf J}_2\\
{\boldsymbol \kappa} 
\end{array}
\right )
&= \sum_{l_{1}} \sum_{l_{23}} \bm{\Theta}^{(1,2, 3)}_{l_{1},l_{23}}
 \left 
 \langle \bU_{l_{23}}^{(2,3)}, 
 \left(\begin{array}{c}
{\bf 1}_{v_1} \otimes {\bf 1}_{v_2} \\
{\bf v}_1\otimes {\bf 1}_{v_2}\\
{\bf 1}_{v_1}\otimes{\bf v}_2\\
\frac12{\bf v}_1^2\otimes {\bf 1}_{v_2}+\frac12{\bf 1}_{v_1}\otimes{\bf v}_2^2
\end{array}
\right )
\right \rangle
\bU_{l_{1}}^{(1)},
\label{eq:rho_j_kappa_2d2v}
\end{align}
where $\bU^{(2,3)}$ is given in \eqref{eq: U34_htd}.

\item {\bf Project.} Then ${\bf f}$ is projected onto subspace  ${\mathcal{N}}$, denoted by ${\bf f}_1 =  {P}_{\mathcal{N}}({\bf f})$ (consistently with the subscript $1$ in the notations) that can preserve the macroscopic densities of ${\bf f}$. With the consideration of the weighted inner product, we need the scaling and rescaling processes, see \cite{guo2022conservative}.

We first construct a set of orthonormal basis of $\mathcal{N}$, denoted by
 $\{\bV_1, \bV_2,\bV_3,\bV_4\}$ in the $(v_1,v_2)$ dimensions from a set of orthonormal basis for $v_1$ and $v_2$ directions as
%
\begin{eqnarray}
\bV_1 &=& \frac{1}{c_1^2} {\bf 1}_{v_1}\otimes {\bf 1}_{v_2},  \quad
\bV_2 = \frac{1}{c_1 c_2} {\bf v}_1\otimes {\bf 1}_{v_2}, \quad
\bV_3 =  \frac{1}{c_1 c_2} {\bf 1}_{v_1}\otimes {\bf v}_2, \nonumber\\
\bV_4 &=& \frac{1}{\sqrt{2}}\left(\frac{1}{c_1c_3}\left(({\bf v}_1^2-c{\bf 1}_{v_1})\right) \otimes ({\bf 1}_{v_2})+\frac{1}{c_1c_3}({\bf 1}_{v_1})\otimes\left(({\bf v}_2^2-c {\bf 1}_{v_2})\right)\right), 
\label{eq: V1_4}
\end{eqnarray} 
with constant $c= \frac{\langle \mathbf{1}_{v_1},\bv_1^2 \rangle_{{\bf w}^{(1)}} }{\langle {\mathbf{1}_{v_1},\mathbf{1}_{v_1}}\rangle_{{\bf w}^{(1)}}}$ for orthogonalization of the basis. $c_l$, $l=1, 2, 3$ are normalization constants for the corresponding basis of ${\bf 1}_{v_1}$, ${\bf v}_1$ and ${\bf v}_1^2-c  {\bf 1}_{v_1}$, where we have assumed the same weight function and discretization in $v_1$ and $v_2$ directions for simplicity. 

The rescaled orthogonal projection is given by 
\begin{equation}
\label{eq:htd_f1_nested}
{\bf f}_1 =  \sum_{l=1}^{4}
(\bU_1^{(1)})_{l}\otimes (\bU_1^{(2,3)})_{l},
\end{equation}
where, 
 $(\bU_1^{(1)})_{l}$, $l=1, \cdots, 4$, are given as
\beq
\label{eq: U12_htd_b}
(\bU_1^{(1)})_{1} =\frac{1}{c_1^2} {\boldsymbol\rho}, \quad
(\bU_1^{(1})_{2} =\frac{1}{c_1 c_2} {\bf J}_1, \quad
(\bU_1^{(1)})_{3} =\frac{1}{c_1 c_2}  {\bf J}_2, \quad
(\bU_1^{(1)})_{4} =\frac{\sqrt{2}}{c_1c_3}({\boldsymbol\kappa} - c\boldsymbol\rho),
\eeq
with  $\bm{\rho}$, $\bJ_1$,  $\bJ_2$, and $\bm{\kappa}$ being the discrete macroscopic charge, current and kinetic energy densities of ${\bf f}$ given in \eqref{eq:rho_j_kappa_2d2v}, and 
\begin{eqnarray}
(\bU_1^{(2,3)})_1 &=& \frac{1}{c_1^2} ({\bf w}^{(1)}\star{\bf 1}_{v_1})\otimes ({\bf w}^{(2)}\star{\bf 1}_{v_2}),  \nonumber\\
(\bU_1^{(2,3)})_2 &=& \frac{1}{c_1 c_2} ({\bf w}^{(1)}\star{\bf v}_1)\otimes ({\bf w}^{(2)}\star{\bf 1}_{v_2}), \nonumber\\
(\bU_1^{(2,3)})_3 &=& \frac{1}{c_1 c_2} ({\bf w}^{(1)}\star{\bf 1}_{v_1})\otimes ({\bf w}^{(2)}\star{\bf v}_2), \nonumber\\
(\bU_1^{(2,3)})_4 &=& \frac{1}{\sqrt{2}}\left(\frac{1}{c_1c_3}\left({\bf w}^{(1)}\star({\bf v}_1^2-c {\bf 1}_{v_1})\right) \otimes ({\bf w}^{(2)}\star{\bf 1}_{v_2})+\frac{1}{c_1c_3}({\bf w}^{(1)}\star{\bf 1}_{v_1})\otimes\left({\bf w}^{(2)}\star({\bf v}_2^2-c {\bf 1}_{v_2})\right)\right). \nonumber\\
\label{eq: U1_4}
\end{eqnarray} 
Hence, we construct the three frame vectors for node $(2)$ as
\begin{eqnarray}
(\bU_1^{(2)})_1 = \frac{1}{c_1} {\bf w}^{(1)}\star{\bf 1}_{v_1} , \quad
(\bU_1^{(2)})_2 = \frac{1}{c_2}{\bf w}^{(1)}\star{\bf v}_1 , \quad
(\bU_1^{(2)})_3 = \frac{1}{c_3}{\bf w}^{(1)}\star({\bf v}_1^2-c {\bf 1}_{v_1}).
\label{eq: U3U4}
\end{eqnarray}
We have the same three frame vectors for the node $(3)$ but for $v_2$, again assuming that the weight function and discretization in $v_2$ is the same as $v_1$, 
\begin{eqnarray}
(\bU_1^{(3)})_1 = \frac{1}{c_1} {\bf w}^{(2)}\star{\bf 1}_{v_2} , \quad
(\bU_1^{(3)})_2 = \frac{1}{c_2}{\bf w}^{(2)}\star{\bf v}_2 , \quad
(\bU_1^{(3)})_3 = \frac{1}{c_3}{\bf w}^{(2)}\star({\bf v}_2^2-c {\bf 1}_{v_2}).
\label{eq: U3U4_b}
\end{eqnarray} 
The transfer tensor $\bm{\Theta}_1^{(2,3)}$ is a third order tensor of size $3\times 3 \times4$. It has zero elements, except the following specification for $(\bm{\Theta}_1^{(2,3)})_{l_{3},l_{4}, l_{34}}$
\begin{equation}
\label{eq:B34f}
(\bm{\Theta}_1^{(2,3)})_{1,1,1} =(\bm{\Theta}_1^{(2,3)})_{2,1,2}=(\bm{\Theta}_1^{(2,3)})_{1,2,3}= 1, \quad
(\bm{\Theta}_1^{(2,3)})_{3,1,4} = (\bm{\Theta}_1^{(2,3)})_{1,3,4} = \frac1{\sqrt{2}}.
\end{equation} 


\item {\bf Truncate in $l^2_{\bf w}$.} We perform a weighted hierarchical HOSVD  truncation to the remainder ${\bf f}_2 = \mathbf{f} -  {\bf f}_1$, followed by a projection operator $(I-P_{\mathcal{N}})$ to ensure zero charge, current, and kinetic energy densities. That is, we  compute $(I-P_{\mathcal{N}}) (\sqrt{\bf w}\star \mathcal{T}_\varepsilon (\frac{1}{\sqrt{\bf w}}\star{\bf f}_2))$, where $\mathcal{T}_\varepsilon $ denotes the standard hierarchical HOSVD  truncation procedure with threshold $\varepsilon$.

\item {\bf Update.} We obtain the low rank truncation of ${\bf f}$ with local mass, momentum and energy conservation, denoted as  
\beq
T_c({\bf f}) = {\bf f}_1 + (I-P_{\mathcal{N}}) (\sqrt{\bf w}\star \mathcal{T}_\varepsilon (\frac{1}{\sqrt{\bf w}}\star{\bf f}_2)). 
\label{eq: Tc}
\eeq
We call the proposed truncation \eqref{eq: Tc} the conservative truncation, as $T_c({\bf f})$ exactly preserves the charge, current and kinetic energy densities of ${\bf f}$. 
\end{enumerate}

The proposed truncation algorithm removes the redundancy in $\mathbf{f}$ and at the same time preserves macroscopic charge, current and kinetic energy densities as proved in \cite{guo2022conservative}.

  \subsection{Local macroscopic conservation achieved by a conservative kinetic flux vector splitting scheme for macroscopic equations}
  \label{sec: LoMac1D1V}

In this subsection, we discuss the technique to attain the LaMoC property under the low rank tensor framework, as an extension of our previous work \cite{guo2022local} to the VM system. Note that
the conservative truncation introduced above indicates that  ${\bf f}_1$ is uniquely determined by macroscopic $\bm{\rho}$, $\bJ_1$,  $\bJ_2$ and $\boldsymbol\kappa$, see \eqref{eq: U12_htd_b}, and the truncated low rank solution $T_c({\bf f})$ preserves these macroscopic densities. On the other hand, it has been well known that numerical methods for the system of conservation laws \eqref{eq:mass1d}-\eqref{eq:ener1d}, if being written in the flux-difference form, can locally conserve the macroscopic quantities. 

To achieve the LoMaC property, the proposed algorithm updates the macroscopic densities $\bm{\rho}$, $\bJ_1$,  $\bJ_2$ and $\boldsymbol\kappa$ by  a high order conservative finite difference discretization of macroscopic system \eqref{eq:mass1d}-\eqref{eq:ener1d}. Since the Vlasov solution ${\bf f}$ is known, the numerical fluxes can be computed by taking the upwind components and performing integration of the flux functions in velocity directions as in KFVS \cite{mandal1994kinetic, xu1995gas}. Once these macroscopic quantities are computed, they are plugged into \eqref{eq: U12_htd_b} to construct a new ${\bf f}_1^{M}$ 
which then replaces ${\bf f}_1$ in \eqref{eq: Tc} to update the low rank Vlasov solution. By design, the remainder part ${\bf f}_2 ={\bf f}-{\bf f}_1$ stays the same with zero macroscopic $\boldsymbol\rho$, $\bJ_1$, $\bJ_2$ and $\boldsymbol\kappa$. In other words, we perform a correction step on the first few moments of ${\bf f}$, by taking advantage of a conservative KFVS scheme for macroscopic equations, to ensure local macroscopic conservation. 

Below we describe the conservative update of macroscopic variables, denoted as $\boldsymbol\rho^{n+1, M}$, $\bJ_1^{n+1, M}$, $\bJ_2^{n+1, M}$, $\boldsymbol\kappa^{n+1, M}$. Let $U \doteq (\rho, P_1,P_2, {e})^\top$,  $F \doteq (J_1,\sigma_1, \sigma_2,   Q)^\top$, then the macroscopic system \eqref{eq:mass1d}-\eqref{eq:ener1d} is written in the following compact form
\beq
\label{eq:U}
U_t + F_x = S,
\eeq
where 
\begin{equation}
\label{eq:source}
S=\left(\begin{array}{c}
0\\
\frac12\frac{\partial}{\partial x_1}\left(|E_1|^2 - |E_2|^2 - |B_3|^2\right) + \rho_iE_1\\
\frac{\partial}{\partial x_1}\left(E_1E_2\right) + \rho_iE_2 \\
\frac{\partial}{\partial x_1}\left(E_2B_3\right)
\end{array}
\right ).
\end{equation}
Assuming the same spatial grid \eqref{eq: x_grid},  the algorithm with the high order upwind finite difference spatial discretization coupled with the second order SSP multi-step time integrator for system \eqref{eq:U} writes
\begin{equation}
\label{eq:Uupdate}
U_j^{n+1} = \frac14U^{n-2}_j + \frac34U^{n}_j + \frac32\Delta t\left(-\frac{1}{h_{x_1}} \left( \hat{F}^n_{j+\frac12} -\hat{F}^n_{j-\frac12}\right) + S_j^n\right),
\end{equation}
where $U^{n}_j = (\rho_j^n, P^n_{1,j}, P^n_{2,j},e^{n}_j)^\top$, $j=1,\ldots,N_{x_1}$. 
The numerical fluxes are uniquely defined at cell interfaces and are given by the following upwind splitting 
\begin{equation}
\hat{F}^n_{j+\frac12} = \hat{F}^{n, +}_{j+\frac12} + \hat{F}^{n, -}_{j+\frac12}, \quad j=1,\ldots,N_{x_1}.
\end{equation}
To obtain $\hat{F}^{n, \pm}_{j+\frac12}$ with high order spatial accuracy in an upwind fashion, assuming the Vlasov solution ${\bf f}^n$ in the low rank HT format \eqref{eq:htd_f_nested_0}-\eqref{eq: U34_htd},  we first compute ${\bf F}^{n, +}$ and ${\bf F}^{n, -}\in \mathbb{R}^{N_{x_1}}$
 \begin{align}
{\bf F}^{n, +}
 &= \sum_{l_{1}} \sum_{l_{23}} \bm{\Theta}^{n,(1,2, 3)}_{l_{1},l_{23}}
 \left 
 \langle \bU_{l_{23}}^{n,(2,3)}, 
 \left(\begin{array}{c}
\bv_1^+\otimes {\bf 1}_{v_2} \\
(\bv_1^+)^2\otimes {\bf 1}_{v_2}\\
\bv_1^+ \otimes \bv_2\\
\frac12 (\bv^+_1)^3\otimes {\bf 1}_{v_2} + \frac12 \bv^+_1\otimes (\bv_2)^2
\end{array}
\right )
\right \rangle
\ \bU_{l_{1}}^{n,(1)}, \\
{\bf F}^{n, -}
 &= \sum_{l_{1}} \sum_{l_{23}} \bm{\Theta}^{n,(1,2, 3)}_{l_{1},l_{23}}
 \left 
 \langle \bU_{l_{23}}^{n,(2,3)}, 
 \left(\begin{array}{c}
\bv_1^-\otimes {\bf 1}_{v_2} \\
(\bv_1^-)^2\otimes {\bf 1}_{v_2}\\
\bv_1^- \otimes \bv_2\\
\frac12 (\bv^-_1)^3\otimes {\bf 1}_{v_2} + \frac12 \bv^-_1\otimes (\bv_2)^2
\end{array}
\right )
\right \rangle
\ \bU_{l_{1}}^{n,(1)},
\label{eq:Fpm}
\end{align}
where $\bv^+_1 = \max(\bv_1, 0)$, $\bv^-_1 = \min(\bv_1, 0)$ and the inner product $\langle \cdot, \cdot \rangle$ is in the sense of \eqref{eq: vm_inner_prod_2_d}. ${\bf F}^{n, \pm}$ can be computed by the standard tensor-vector contraction.
Let $F^{n, \pm}_j ={\bf F}^{n, \pm}(j)$, the upwind fluxes $\hat{F}^{n, \pm}_{j+\frac12}$ are reconstructed from ${\bf F}^{n, \pm}(:)$ in the following way using the fifth order upwind stencils \cite{shu2009high}, 
\begin{align*}
\hat{F}^{n, -}_{j+\frac12} &= -\frac{1}{20}F_{j-1}^{n, -} + \frac{9}{20}F_{j}^{n, -} + \frac{47}{60}F_{j+1}^{n, -} - \frac{13}{60}F_{j+2}^{n, -}+ \frac{1}{30}F_{j+3}^{n, -},\\
 \hat{F}^{n, +}_{j+\frac12} &= \frac{1}{30}F_{j-2}^{n, +} -\frac{13}{60}F_{j-1}^{n, +} + \frac{47}{60}F_{j}^{n, +} +  \frac{9}{20}F_{j+1}^{n, +} -\frac{1}{20}F_{j+2}^{n, +}.
\end{align*}
Furthermore, we employ the spectral method to approximate the derivatives in source term $S$ \eqref{eq:source}. Then we let 
$
 U^{n+1}_j = \left(\begin{array}{l}
{\rho}_j^{n+1, M} \\
P_{1,j}^{n+1, M}\\
P_{2,j}^{n+1, M}\\
 {e}_j^{n+1, M}
\end{array}
\right )$ updated from \eqref{eq:Uupdate}, from which we can compute 
\begin{align}
J_{1,j}^{n+1, M} &= {P}_{1,j}^{n+1, M} -E_{2,j}^{n+1}B_{3,j}^{n+1},\label{eq:kinetic_updatej1}\\
J_{2,j}^{n+1, M} &= {P}_{2,j}^{n+1, M} +E_{1,j}^{n+1}B_{3,j}^{n+1}\label{eq:kinetic_updatej2},\\
\label{eq:kinetic_updatek}
{\kappa}_j^{n+1, M} &= {e}_j^{n+1, M} - \frac12\left(|E_{1,j}^{n+1}|^2 + |E_{2,j}^{n+1}|^2+|B_{3,j}^{n+1}|^2\right).
\end{align}
The electric and magnetic fields $\bE^{n+1}_1$, $\bE^{n+1}_2$, $\bB^{n+1}_3$ are computed  via Maxwell's equations. Note that $\bE^{n+1}_1$ can be updated using Ampère's law in \eqref{eq:max11d} or Gauss's law \eqref{eq:gauss1d}. In the simulation, we solve Poisson's equation with charge density $\boldsymbol\rho^{n+1, M}$ to update $\bE^{n+1}_1$, and hence Gauss's law is satisfied. Lastly, we construct ${\bf f}^M_1$ in the same way as ${\bf f}_1$ using \eqref{eq:htd_f1_nested}, except the basis $\bU^{(1)}$ in \eqref{eq: U12_htd_b} is defined using $\boldsymbol\rho^{n+1, M}$, $\bJ_1^{n+1, M}$, $\bJ_2^{n+1, M}$, $\boldsymbol\kappa^{n+1, M}$ updated from the macroscopic system.  Such a replacement can be viewed as a correction step for macroscopic conservation. 
Meanwhile, the treatment for ${\bf f}_2$ in the orthogonal decomposition  stays the same. That is ${\bf f}_2$ is truncated but still preserves zero charge, current and kinetic energy densities after truncation. Last, the LoMaC low rank solution is updated as 
\begin{equation}
\label{eq:lomac}
{\bf f}^{n+1} \doteq T^M_c({\bf f}) = {\bf f}^M_1 +(I-P_\mathcal{N})\left(\sqrt{\bf w}  \star \mathcal{T}_\varepsilon (\frac{1}{\sqrt{\bf w}} \star {\bf f}_2)\right).
\end{equation}

We summarize the newly proposed LoMaC low rank tensor algorithm for the 1D2V VM system in Algorithm \ref{alg: low_rank_lomac_1D2V}, based on a low rank Vlasov solver using  the spectral collocation spatial discretization and second order SSP multi-step temporal discretization, together with a conservative high order finite difference macroscopic solver with KFVS and the same temporal discretization.

\bigskip
 \bigskip
\begin{algorithm}[H]
\label{alg: low_rank_lomac_1D2V}
  \caption{The LoMaC low rank tensor algorithm for the 1D2V VM system.}
    \begin{enumerate}
\item Initialization:
  \begin{enumerate}
  \item Initial distribution function $\mathbf{f}^0$ in a low rank format \eqref{eq:htd_f_nested_0}-\eqref{eq: U34_htd} together with electric and magnetic fields $\bE^0_1$, $\bE^0_2$, $\bB_3^0$ and macroscopic densities $\boldsymbol\rho^{0}$, $\bP_1^{0}$, $\bP_2^{0}$, ${\boldsymbol e}^{0}$.
 \end{enumerate}
\item For each time step evolution from $t^n$ to $t^{n+1}$: update ${\bf f}^{n+1}$ from ${\bf f}^n$ in the low rank HT format.
  \begin{enumerate}
  \item Update ${\bf f}^{n+1,*}$, by adding basis according to the dimension tree $\mathcal{T}$ as shown in Figure \ref{fig:dimtree1}. The procedure is similar to that outlined in \cite{guo2021lowrank}. 
  \item  Compute ${\bf f}_2\doteq {\bf f}^{n+1,*} - {\bf f}_1$ and perform a weighted hierarchical  HOSVD truncation on ${\bf f}_2$ in the low rank 1D2V format \cite{guo2022conservative} to obtain
  $
\left( \sqrt{\bf w}  \star \mathcal{T}_\varepsilon (\frac{1}{\sqrt{\bf w}} \star {\bf f}_2)\right).
  $
Finally, we apply the $(I-P_\mathcal{N})$ operator to $\left( \sqrt{\bf w}  \star \mathcal{T}_\varepsilon (\frac{1}{\sqrt{\bf w}} \star {\bf f}_2)\right)$ to ensure its zero charge, current and kinetic energy densites after truncation.
  \item Compute ${\bf f}^M_1$.
  \begin{enumerate}
  \item Update macroscopic mass, momentum and energy densites $\boldsymbol\rho^{n+1, M}$, $\bP_1^{n+1, M}$, $\bP_2^{n+1, M}$, ${\boldsymbol e}^{n+1, M}$, with KFVS, in a flux-difference form using the same second order SSP multi-step method in Step 2(a). 
  \item Update $\bE_1^{n+1}$, $\bE_2^{n+1}$, $\bB_3^{n+1}$ from Maxwell's equations. 
  \item Compute  $\bJ_1^{n+1, M}$,  $\bJ_2^{n+1, M}$, and $\boldsymbol\kappa^{n+1, M}$ from \eqref{eq:kinetic_updatej1}-\eqref{eq:kinetic_updatek}. 
\item  Construct ${\bf f}^M_1$ from \eqref{eq: U12_htd_b}, but with $\boldsymbol\rho^{n+1, M}$, $\bJ_1^{n+1, M}$, $\bJ_2^{n+1, M}$, $\boldsymbol\kappa^{n+1, M}$. 
\end{enumerate}
  \item Update the compressed low-rank solution via \eqref{eq:lomac}, 
  $${\bf f}^{n+1} = {\bf f}^M_1 +(I-P_\mathcal{N})\left(\sqrt{\bf w}  \star \mathcal{T}_\varepsilon (\frac{1}{\sqrt{\bf w}} \star {\bf f}_2)\right).$$
\end{enumerate}
\end{enumerate}
  \end{algorithm}

The proposed LoMaC low rank update of the VM solution starts with the adding basis step. In particular, the method employs a spectral collocation spatial discretization together with an SSP multi-step time integrator, ensuring high order accuracy in both space and time. By simultaneously updating macroscopic conservation laws with KFVS, together with a weighted orthogonal projection of the low rank solution, the method can preserve the same macroscopic mass, momentum and energy density as the conservation laws. Last, we apply a conservative weighted hierarchical HOSVD truncation to remove redundancy and at the same time retain the conservation properties.  Note that for one step evolution,  macroscopic and kinetic parts are independent except using $\boldsymbol\rho^{n+1, M}$, $\bJ_1^{n+1, M}$, $\bJ_2^{n+1, M}$, $\boldsymbol\kappa^{n+1, M}$ to construct ${\bf f}^M_1$ in \eqref{eq: U12_htd_b} and using $\mathbf{f}^n$ to construct the numerical flux with KFVS. We summarize the LoMaC property of the proposed low rank  method in the following algorithm.

\begin{prop} The proposed LoMaC low rank VM solver locally conserves the macroscopic mass, momentum and energy. 
\end{prop}
\begin{proof} The proof follows directly from the construction of the algorithm.
\end{proof}

%% file: AlgorithmNew.tex
\subsection{A low rank tensor approach for the 2D2V Vlasov-Maxwell system}

We  then consider the  2D2V VM system given as follows:
\begin{equation}
  \frac{\partial f}{\partial t}+v_{1} \frac{\partial f}{\partial x_{1}}+v_{2} \frac{\partial f}{\partial x_{2}}+\left(E_{1}+v_{2} B_{3}\right) \frac{\partial f}{\partial v_{1}}+\left(E_{2}-v_{1} B_{3}\right) \frac{\partial f}{\partial v_{2}}=0,
  \label{vlasov2}
\end{equation}

\begin{equation}
\label{eq:max2d}
     \dfrac{\partial E_{1}}{\partial t}=\dfrac{\partial B_{3}}{\partial x_{2}}-J_1,\quad
     \dfrac{\partial E_{2}}{\partial t}=-\dfrac{\partial B_{3}}{\partial x_{1}}-J_2,\quad
     \dfrac{\partial B_{3}}{\partial t}=\dfrac{\partial E_{1}}{\partial x_{2}}-\dfrac{\partial E_{2}}{\partial x_{1}},
\end{equation}


\begin{equation}
\frac{\partial E_{1}}{\partial x_{1}}+\frac{\partial E_{2}}{\partial x_{2}}=\rho-\rho_i,\label{eq:gauss2d}
\end{equation}
where $\rho_i$ denotes the constant background density of ions.  The associated macroscopic conservation laws become 
\begin{align} 
	&\frac{\partial \rho}{\partial t} + \frac{\partial J_1}{\partial x_1}+ \frac{\partial J_2}{\partial x_2} = 0,\\
 &\frac{\partial P_1}{\partial t}+ \frac{\partial \sigma_{11}}{\partial x_1}+\frac{\partial \sigma_{12}}{\partial x_2} =\frac12\frac{\partial}{\partial x_1}\left(|E_1|^2 - |E_2|^2 - |B_3|^2\right) + \frac{\partial}{\partial x_2}\left(E_1E_2\right)+ \rho_iE_1, \\ 
     &\frac{\partial P_2}{\partial t}+ \frac{\partial \sigma_{21}}{\partial x_1}+\frac{\partial \sigma_{22}}{\partial x_2} = \frac{\partial}{\partial x_1}\left(E_1E_2\right)+\frac12\frac{\partial}{\partial x_2}\left(-|E_1|^2 + |E_2|^2 - |B_3|^2\right)+\rho_iE_2,\\
      & \frac{\partial e}{\partial t} + \frac{\partial Q_1}{\partial x_1} + \frac{\partial Q_2}{\partial x_2} =\frac{\partial\left(E_2B_3\right)}{\partial x_1}- \frac{\partial\left(E_1B_3\right)}{\partial x_2},
\end{align}
where the flux functions are defined as $ \sigma_{11} = \int   v_1^2 f \; dv_1dv_2,\; \sigma_{12}=\sigma_{21} = \int   v_1  v_2 f \; dv_1dv_2,\; \sigma_{22} = \int   v_2^2 f\; dv_1dv_2,\; Q_1 =\frac12\int v_1(v_1^2 +v_2^2 ) f\; dv_1dv_2$, and $Q_2 =\frac12\int v_2(v_1^2 +v_2^2 ) f\; dv_1dv_2$.

We can directly extend the Algorithm 1 for simulating the 2D2V VM  system, while we do not explore the decomposition of $x_1$ and $x_2$ but still employ the same dimension tree as the 1D2V VM system given in Figure \ref{fig:dimtree1}. The main reason is that we need to solve the macroscopic model over the full 2D grid in space to attain conservation conveniently. Further, keeping the full grid discretization for the spatial variables  allows for flexibility in handling complex geometry and general boundary conditions, yet it may also lead to increased computational costs, see our previous work \cite{guo2022local} for more details.  In addition, the proposed algorithm for the 2D2V VM system conserve the macroscopic mass, momentum and energy locally as with the Proposition 2.1.

%% file: numerical.tex
\section{Numerical results}

In this section, we present a collocation of numerical tests to demonstrate the performance of the numerical LoMaC  low rank tensor method for solving the reduced 1D2V and 2D2V VM systems. In the simulations, we employ the exponential filter $\sigma(\eta) = \exp(-\alpha\eta^{2p})$ with $p=4$, $\alpha=35$, see \cite{majda1978fourier,gottlieb2001spectral}, to post-process the basis $\bU^{(1)},\bU^{(2)},\bU^{(3)}$ after each time step update. Such a filter introduces numerical viscosity which helps stabilize the numerical method and keep the hierarchical ranks of the solution tensor low. To decrease the adverse 
 effect of the filter for accuracy, it is only turned on if $r\ge15$ where $r$ is the maximum of the hierarchical ranks.
In addition, we let $N_{v_1}=N_{v_2}=:N_v$. We remark that since the condition $\int E_1 dx_1= \int E_2 dx_1 = 0$ is not enforced, the total momentum may not be conserved exactly at the discrete level, while the momentum densities $P_1$ and $P_2$ are still locally conserved by respecting the conservation laws \eqref{eq:mom11d}-\eqref{eq:mom21d} at the discrete level. Further, Gauss's law is satisfied at the discrete level.

\begin{exa}\label{ex:lan}
We consider the Landau-type problem with the initial condition 
    \[f(x_1,v_1,v_2,t=0) = \frac{1}{2 \pi} e^{-\frac{1}{2} (v_{1}^2+v_{2}^2)} (1+\alpha \cos{(kx_1})),\]
where  $k = 0.4$ and $\alpha = 0.01$. The electric field is initialized by solving Gauss’ law \eqref{eq:gauss1d}, yielding 
   \[ E_1(x_1,t=0) = -\frac{\alpha}{k} \sin{(kx_1)}, \quad E_2(x_1,t=0) = 0,\]  
and the magnetic field at $t=0$ is chosen as
    \[ B_3(x_1,t=0) = \frac{\alpha}{k} \sin{(kx_1)},\]
The spatial domain is $\Omega_{x_1}= [0, 2\pi/k]$ and the velocity domain is chosen as $\Omega_{v_1} \times \Omega_{v_2}=[-5, 5]^2$. We first set $\varepsilon = 10^{-5}$ for truncation. In Figure \ref{fig:1}, we present the time histories of electric and magnetic energy, as well as the time histories of the relative deviation of total mass and energy without considering conservative properties. This conserves total mass and energy up to the truncation, $10^{-5}$ caused by the absence of conservation properties in Algorithm 1.  In Figure \ref{fig:2}, we report the time histories of the electric and magnetic energy,  numerical ranks and the time histories of the relative deviation of total mass and energy of the proposed LoMaC low rank solutions for two sets of meshes $N_x\times N_v^2 = 32 \times 128^2$ and $N_x\times N_v^2 = 64 \times 256^2$.  The results agree with those reported in \cite{einkemmer2020low} by the dynamical low rank method that the electric energy decays and exhibits oscillations over time. It is observed that the LoMaC method conserves the total mass and energy up to the machine precision regardless of the mesh size. Furthermore, we observe the hierarchical  ranks of the solution tensor gradually grow in time to capture the dynamics  especially in dimension $v_1$. We measured the algorithm’s execution time in terms of CPU usage. The CPU cost is $830s$ and $1978s$ for the two meshes with a serial implementation, which only approximately doubles with mesh refinement, indicating that the curse of dimensionality is alleviated for solving this problem.    
Then, we  choose a larger truncation threshold $\varepsilon = 10^{-4}$, and the results are summarized in Figure \ref{fig:3}. The time histories of the electric and magnetic energy qualitatively agree with the result computed using $\varepsilon = 10^{-5}$. Meanwhile, the total mass and energy are conserved up to the machine precision as expected, as the LoMaC property of the proposed method is independent of the truncation threshold. Furthermore, with a lager threshold, the hierarchical  ranks become smaller, and computational cost is reduced accordingly, at the expense of some accuracy.

\begin{figure}
			\centering		
    \subfigure[]{\includegraphics[height=40mm]{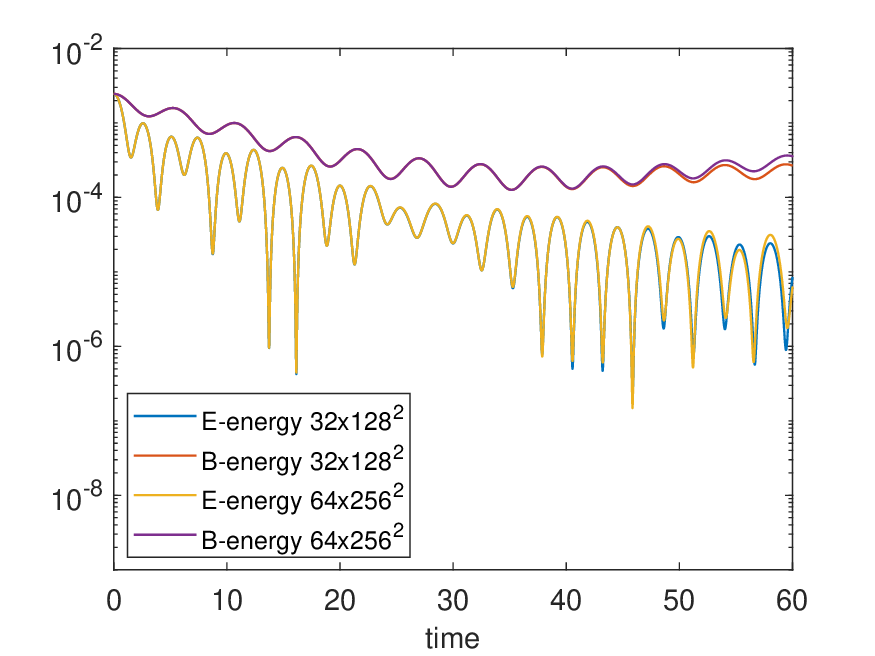}}
	 \subfigure[]{\includegraphics[height=40mm]{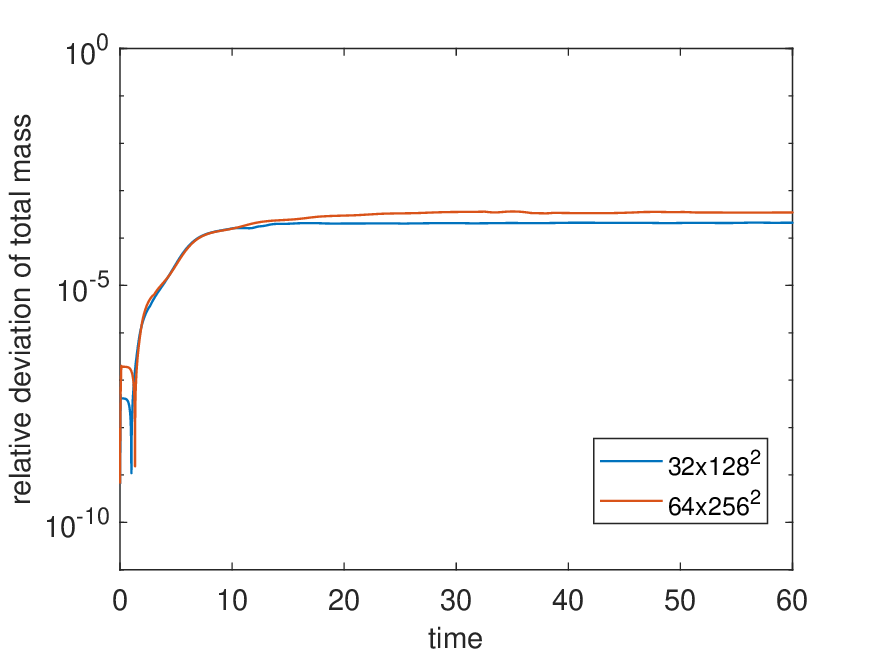}}
     \subfigure[]{\includegraphics[height=40mm]{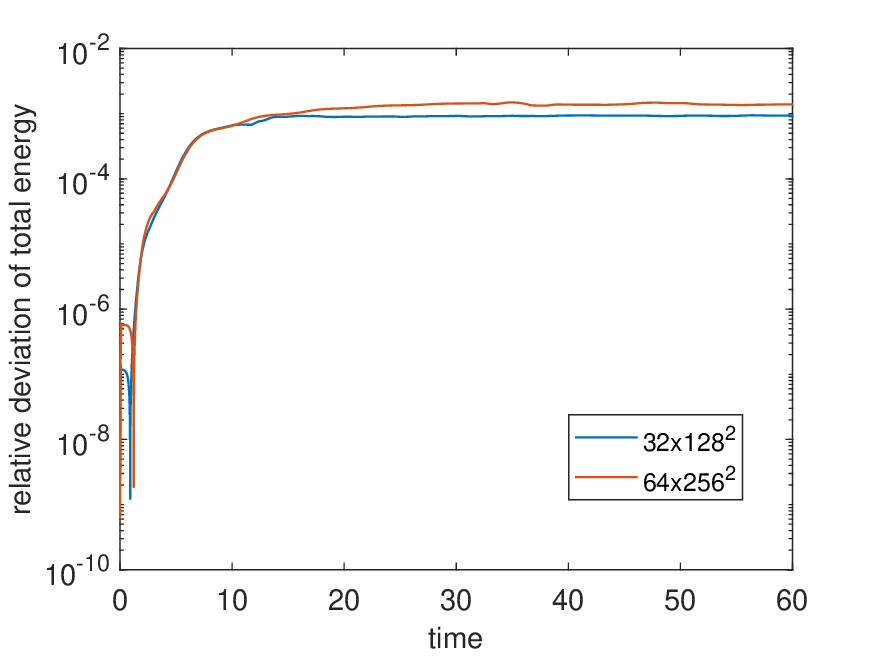}}
    
	\caption{ Example \ref{ex:lan}. The time evolution of the electric and magnetic energy without the conservation (a), relative
deviation of the total mass (d), and total energy  (e). $T=60$ , $ \varepsilon = 10^{-5}$.}
	\label{fig:1}
	\end{figure}	

\begin{figure}
			\centering	
     \subfigure[]{\includegraphics[height=40mm]{Landu_10e5_graph.eps}}
     \subfigure[$N_x\times N_v^2=32 \times 128^2$]{\includegraphics[height=40mm]{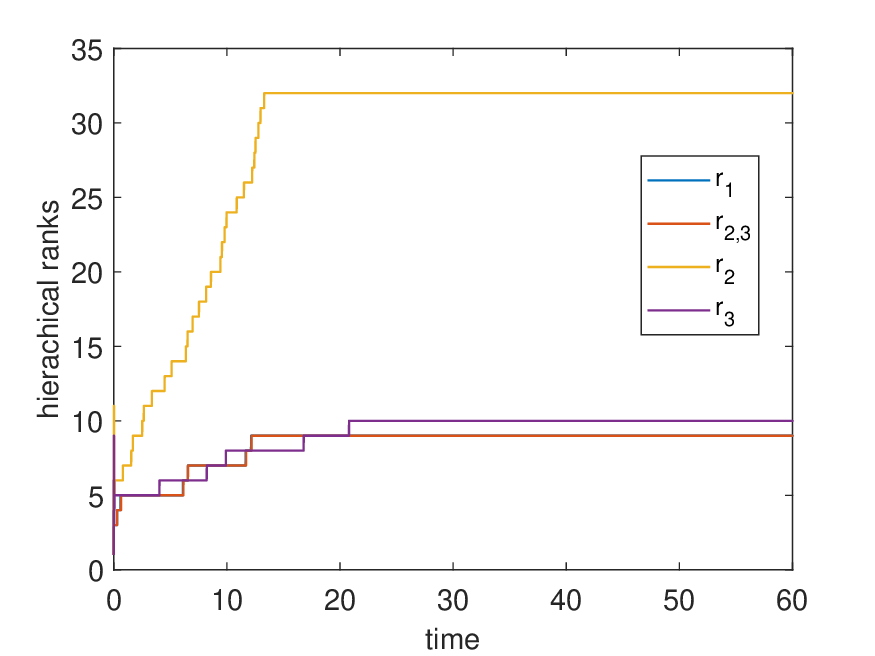}}
      \subfigure[$N_x\times N_v^2=64 \times 256^2$]{\includegraphics[height=40mm]{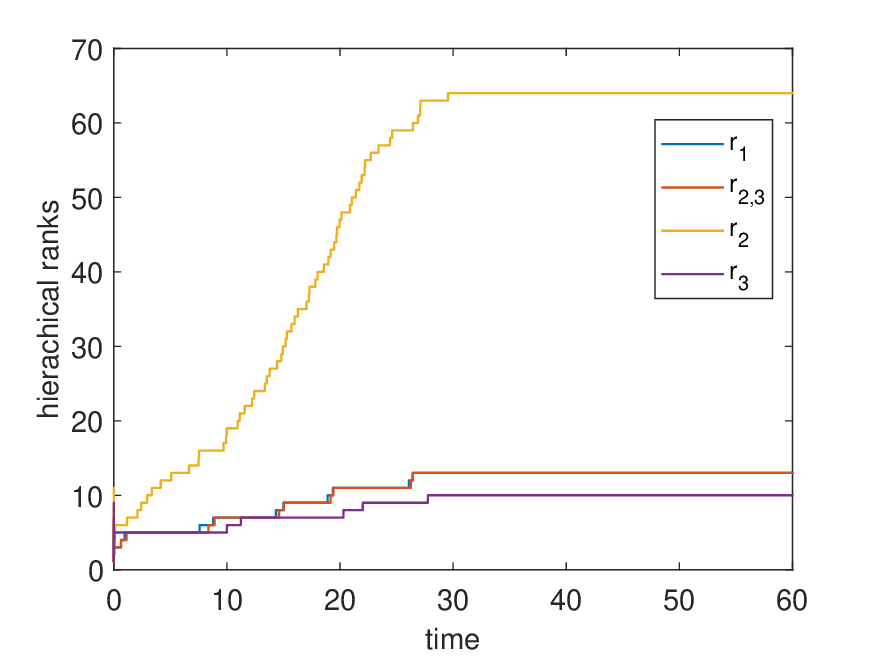}}    
      \subfigure[]{\includegraphics[height=40mm]{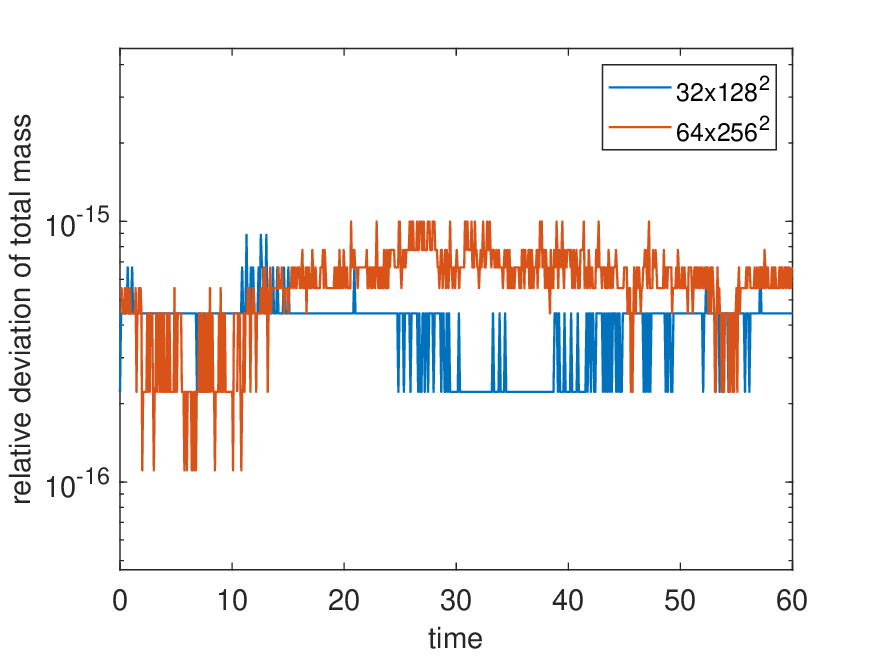}}
     \subfigure[]{\includegraphics[height=40mm]{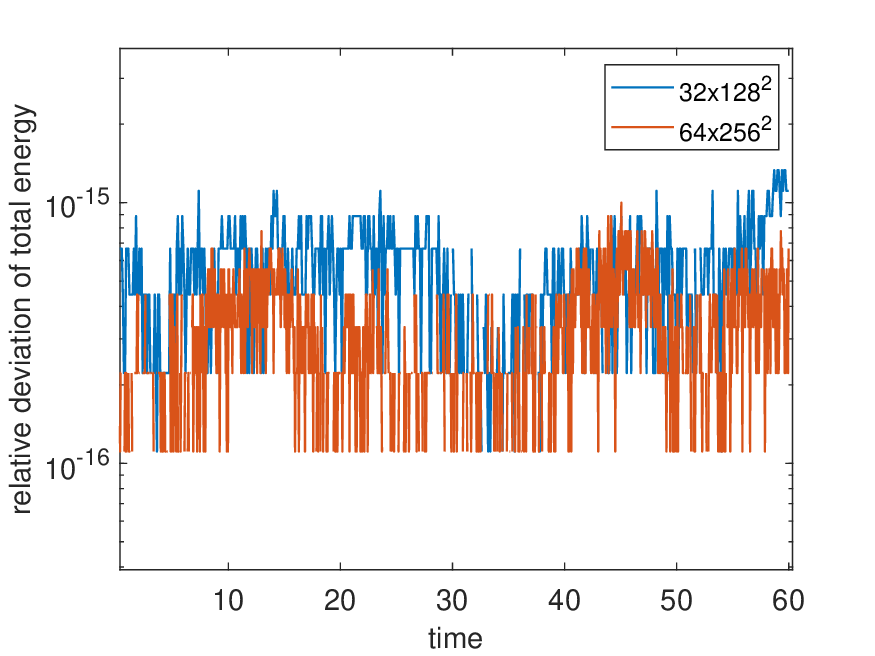}}

	\caption{Example \ref{ex:lan}. The time evolution of the electric and magnetic energy (a), the hierarchical  ranks of the
numerical solution tensor (b, c), relative
deviation of the total mass (d), and total energy  (e). $T=60$. $ \varepsilon = 10^{-5}$. In (b) and (c) the curve of $r_{2,3}$ is close to that of $r_1$.}
	\label{fig:2}
	\end{figure}

\begin{figure}
			\centering	
     \subfigure[]{\includegraphics[height=40mm]{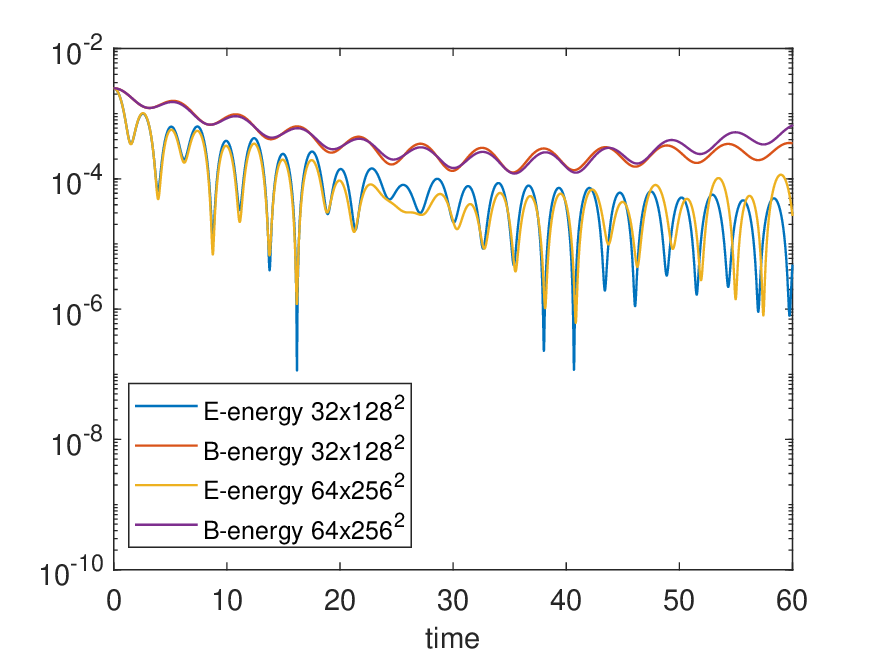}}
       \subfigure[$32 \times 128^2$]
     {\includegraphics[height=40mm]{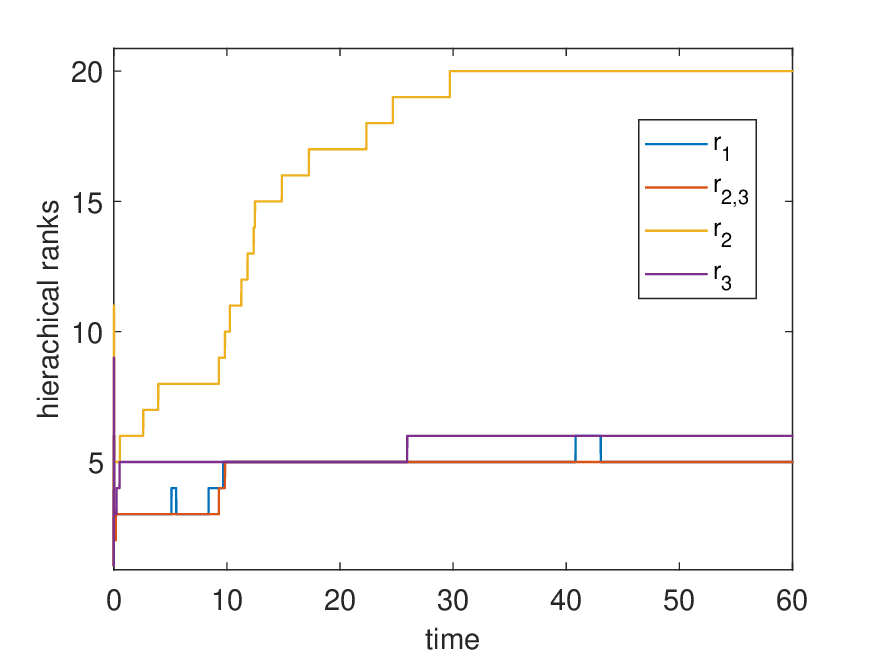}}
      \subfigure[$64 \times 256^2$]{\includegraphics[height=40mm]{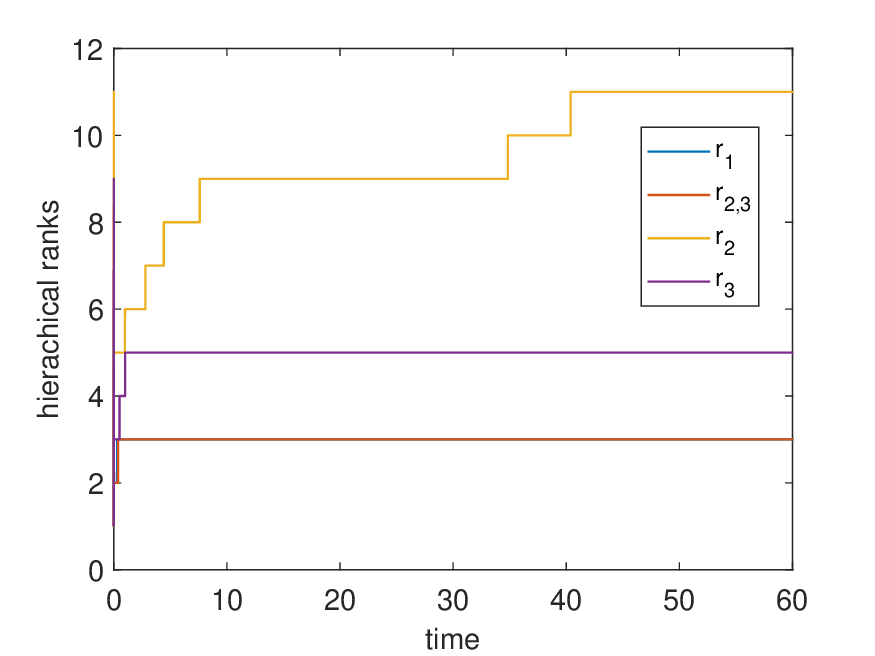}}
     \subfigure[]{\includegraphics[height=40mm]{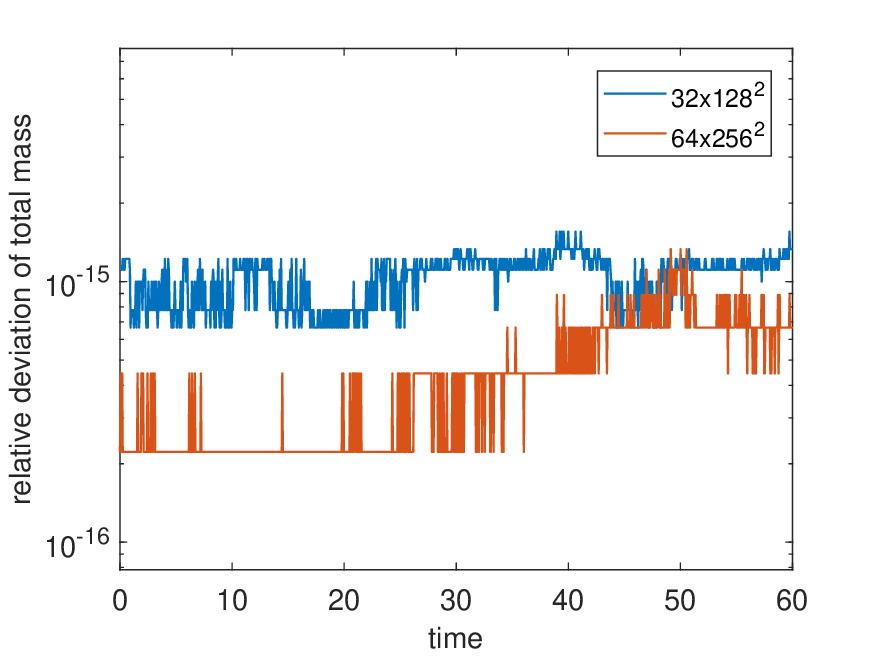}}     
     \subfigure[]{\includegraphics[height=40mm]{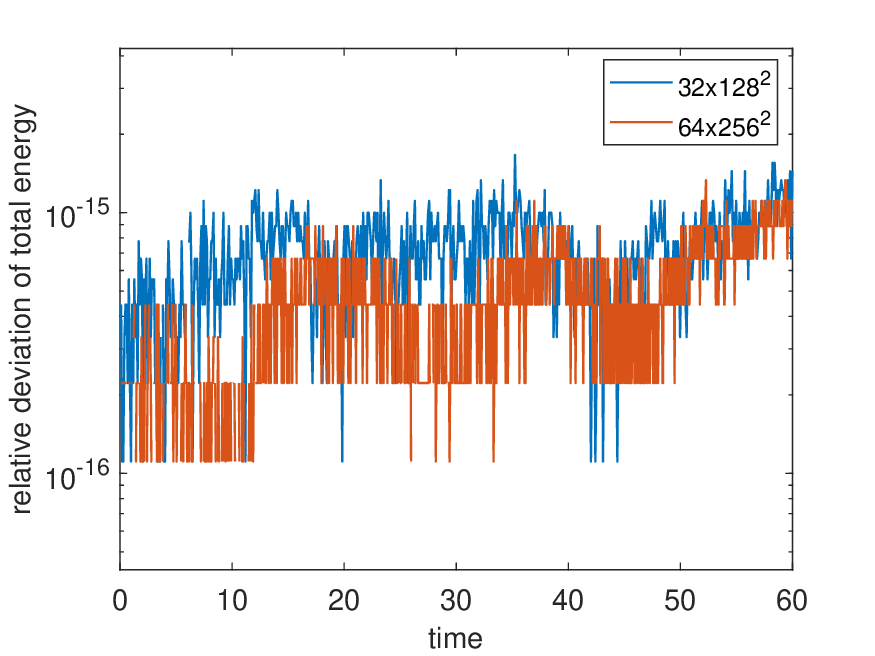}}
	\caption{Example \ref{ex:lan}. The time evolution of the electric and magnetic energy (a), the hierarchical  ranks of the numerical solution tensor (b, c), relative
deviation of the total mass (d), and total energy  (e). $T=60$. $ \varepsilon = 10^{-4}$. In (b) and (c) the curve of $r_{2,3}$ is close to that of $r_1$.}
	\label{fig:3}
	\end{figure}

\end{exa}		
\begin{exa}\label{ex:two}
We consider the two-stream instability problem with initial condition
    \[f(x_1,v_1,v_2,t=0) = \frac{1}{2 \pi \beta} e^{-v_{2}^2/\beta  (e^{-(v_{1}-0.2)^2/\beta} +e^{-(v_{1}+0.2)^2/\beta})},\]
where  $\beta = 2 \cdot10^{-3}$. The spatial domain  is $\Omega_{x_1}= [0, 2\pi]$, and the computational velocity domain is set to be $\Omega_{v_1} \times \Omega_{v_2}=[-0.4, 0.4]^2$.  The electric field is initialized to zero, and the initial magnetic field is chosen as   
\[ B_3(x_1,t=0) = \alpha \sin{(x_1)},\]
with $\alpha = 10^{-3}$, and such perturbation in $B_3$ creates the instability.   We set $\varepsilon = 10^{-5}$ for truncation. In Figure \ref{fig:4}, we depict the time histories of electric, magnetic, and kinetic energy, along with the time histories of the relative deviation of total mass and energy before applying the conservation part. Similar to the previous example, this preserves conservation up to the truncation. In Figure \ref{fig:5}, we plot the time histories of the electric, magnetic, and kinetic energy, hierarchical ranks of the solution tensors, and the time histories of the relative deviation of total mass and energy are displayed for two sets of meshes $N_x\times N_v^2 = 32 \times 64^2$ and $64\times 128^2$. It is observed that the electric energy oscillates in the beginning, then starts to grow exponentially, and eventually saturates.
 The results are consistent with those reported in \cite{einkemmer2020low}. Furthermore, the LoMaC low rank method  conserves the total mass and energy up to the machine precision similar to the previous example.

\begin{figure}
			\centering	
     \subfigure[]{\includegraphics[height=40mm]{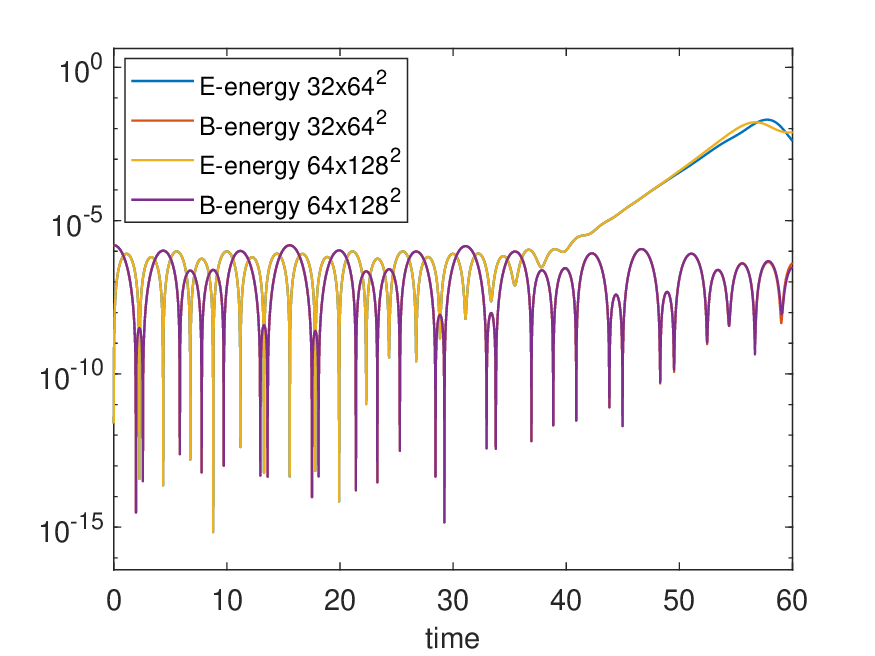}}
	\subfigure[]{\includegraphics[height=40mm]{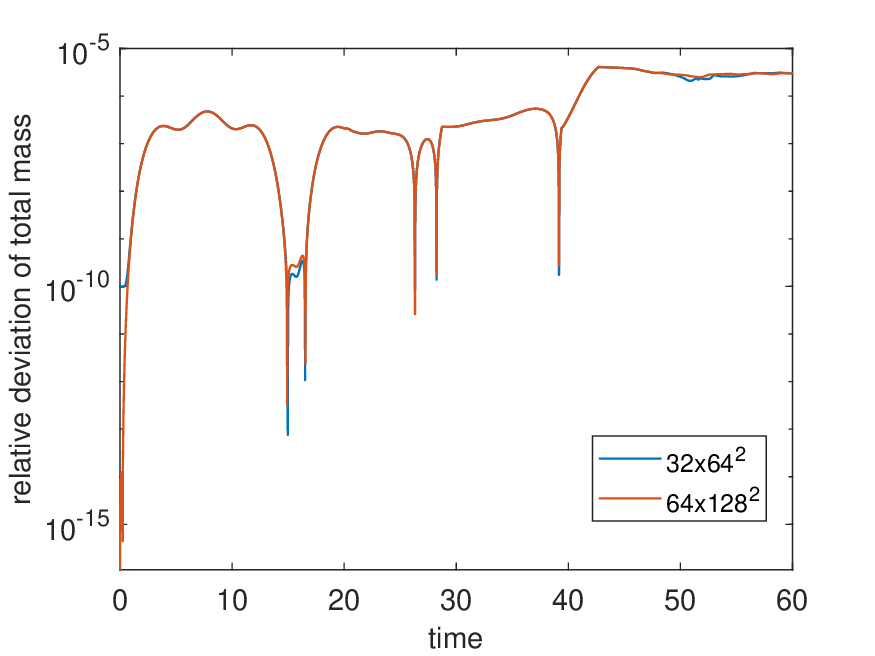}}
     \subfigure[]{\includegraphics[height=40mm]{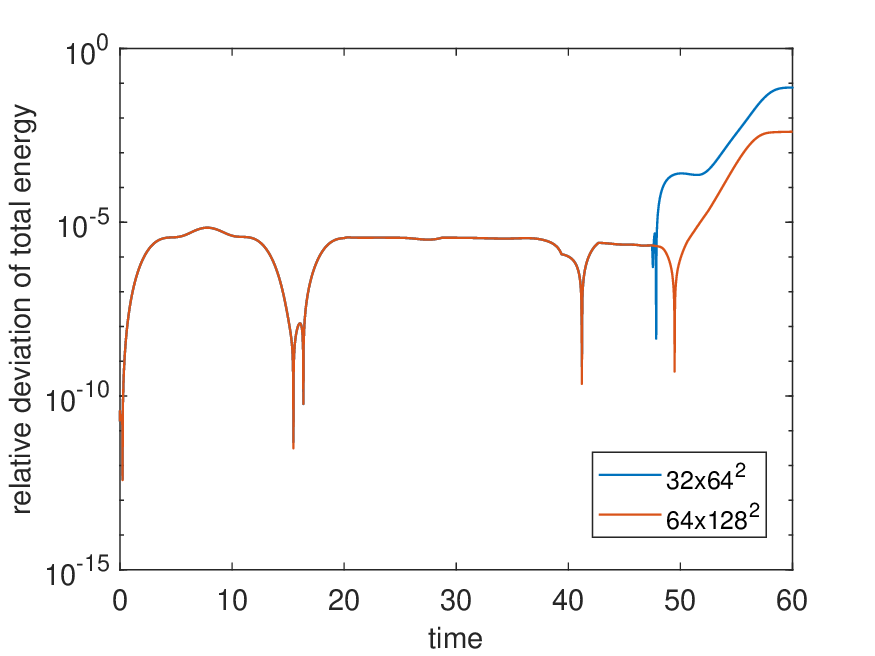}}
     
	\caption{ Example \ref{ex:two}. The time evolution of the electric and magnetic energy without the conservation (a), relative
deviation of total energy (d), and total mass (e). $T=60$ , $ \varepsilon = 10^{-5}$.}
	\label{fig:4}
	\end{figure}

\begin{figure}
			\centering	
     \subfigure[]{\includegraphics[height=40mm]{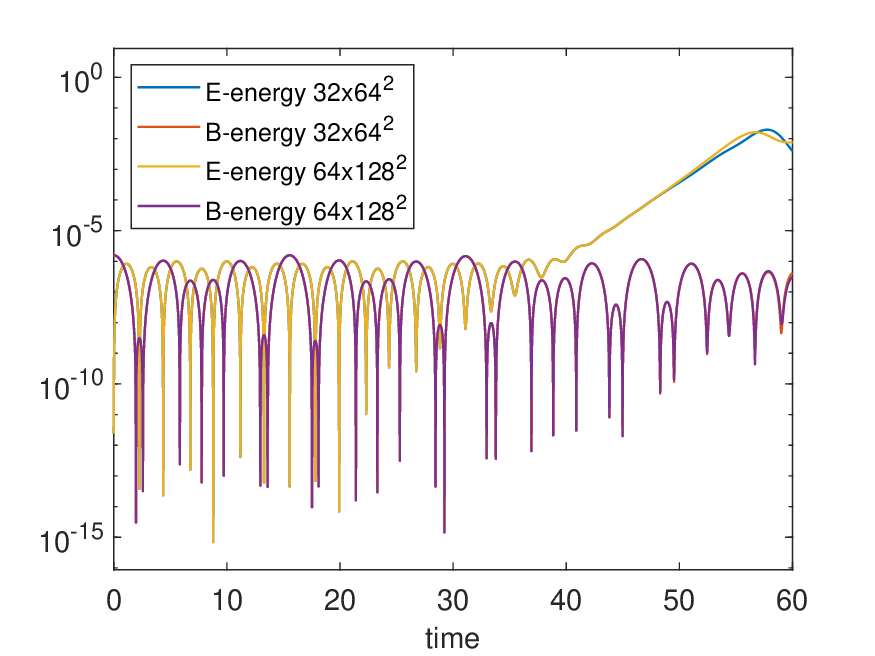}}
       \subfigure[$32 \times 64^2$]
     {\includegraphics[height=40mm]{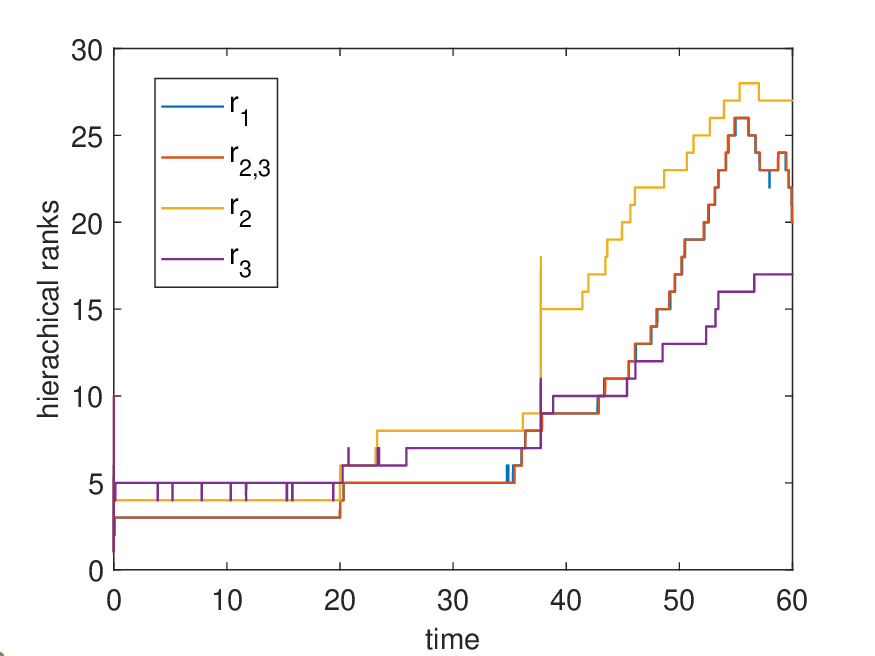}}
      \subfigure[$64 \times 128^2$]{\includegraphics[height=40mm]{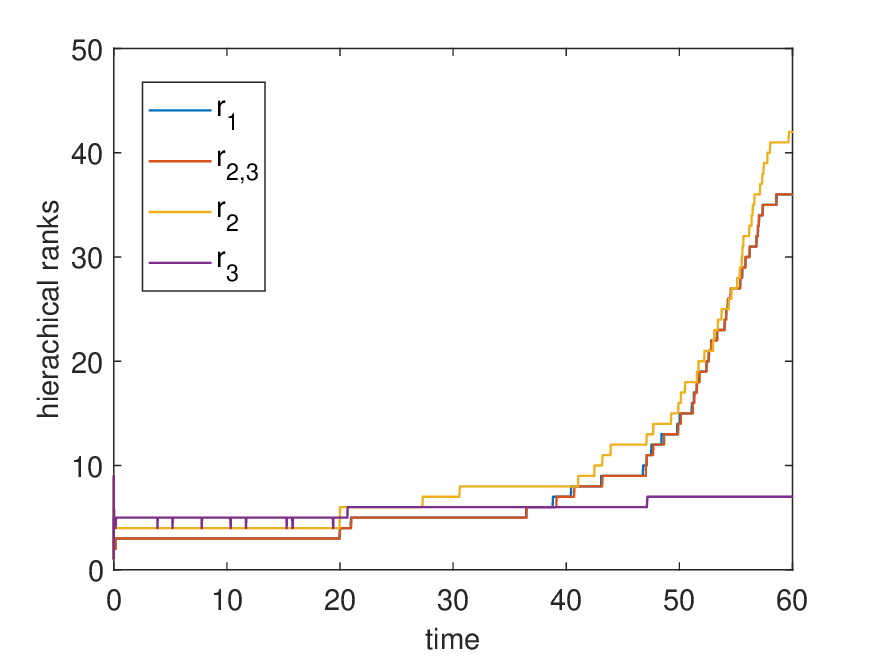}}
     \subfigure[]{\includegraphics[height=40mm]{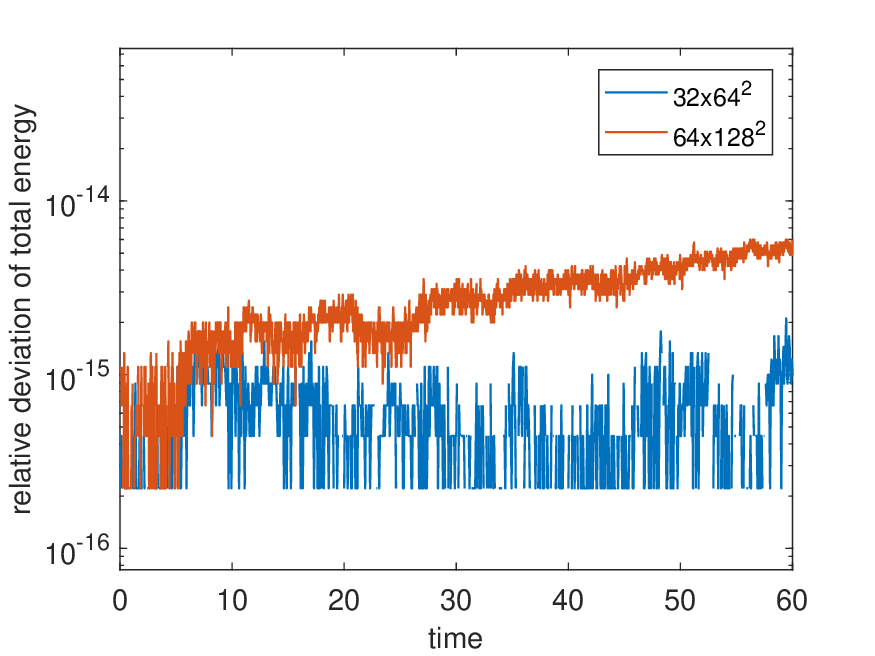}}
     \subfigure[]{\includegraphics[height=40mm]{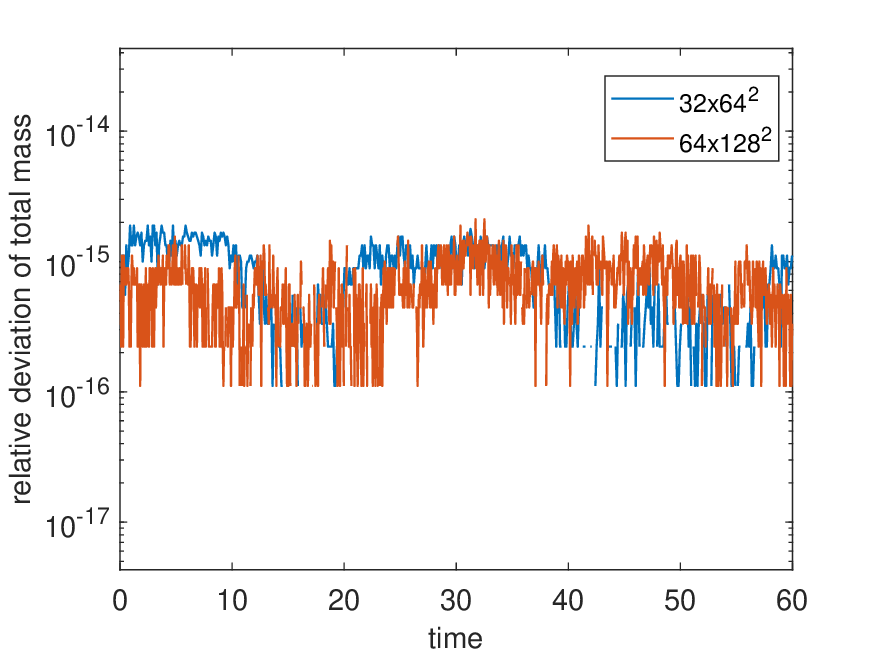}}
	\caption{Example \ref{ex:two}. The time evolution of the electric and magnetic energy for conservative method (a), the rank of the
numerical solutions (b,c), relative
deviation of total energy (d), and total mass (e). $T=60$ , $ \varepsilon = 10^{-5}$. In (b) and (c) the curve of $r_{2,3}$ is close to that of $r_1$.}
	\label{fig:5}
	\end{figure}

\end{exa}

\begin{exa}\label{ex:bump}  We consider the bump-on-tail instability problem for which the initial condition is given by
    \[f(x_1,v_1,v_2,t=0) = \frac{1}{\sqrt{2} \pi} \left(\alpha e^{-v_{1}^2/2} + \beta e^{-2(v_{1}-4.5)^2} (1+\gamma \cos{(kx_1})\right)e^{-v_2^2}, \]
where $\alpha = 9/10$, $\beta = 2/10$, $\gamma = 0.03$  and $k = 0.1$. The computational domain is set to be $\Omega_{x_1} \times \Omega_{v}=[0, 20 \pi] \times [-9, 9]^2$.   The electric field is initialized according to Gauss’s law, and the initial magnetic field is set to be zero.  
In the simulation,  we let $\varepsilon = 10^{-4}$ for truncation. In Figure \ref{fig:8} (a), we plot the time evolution of the electric and magnetic energy. It is observed that, $\bE_2$ and $\bB_3$ stays around the machine precision, and the problem indeed reduces to the 1D1V VP system. In Figure  \ref{fig:8} (b, c), we can see that  the rank $r_3$ corresponding to $v_2$ direction stays close to 1 over time, and hence the proposed low rank tensor method can effective extract the intrinsic 2D structures from the 3D problem. Again, it is observed that  the relative deviation of the total mass and energy is on the scale of the machine precision over time. In Figure \ref{fig:7}, we report 2D cuts of contour plots of the solution at $v_2=0.5$ at several instances of time with mesh $N_x\times N_v^2 = 64 \times 128^2$. It is observed  the numerical results match with those reported in \cite{einkemmer2020low}. 
%

\begin{figure}
			\centering
	  \subfigure[]{\includegraphics[height=40mm]{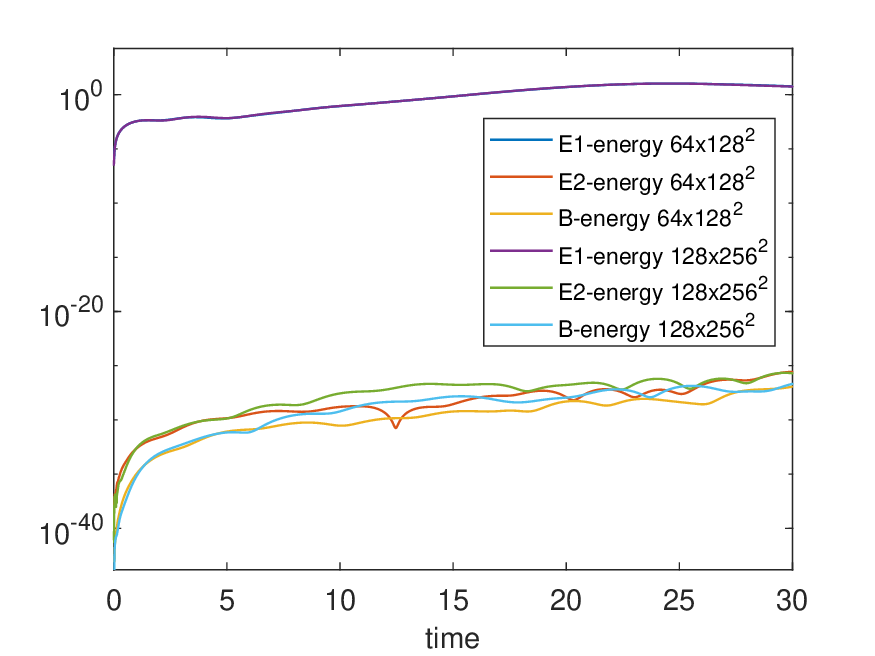}}
            \subfigure[$N_x\times N_v^2=64 \times 128^2$]{\includegraphics[height=40mm]{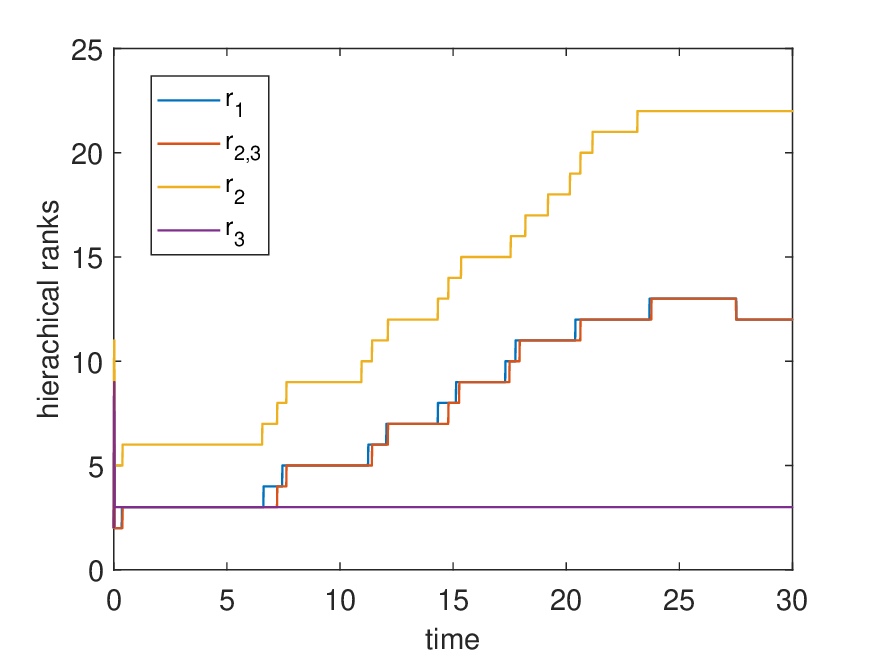}}
            \subfigure[$N_x\times N_v^2=128 \times 256^2$]{\includegraphics[height=40mm]{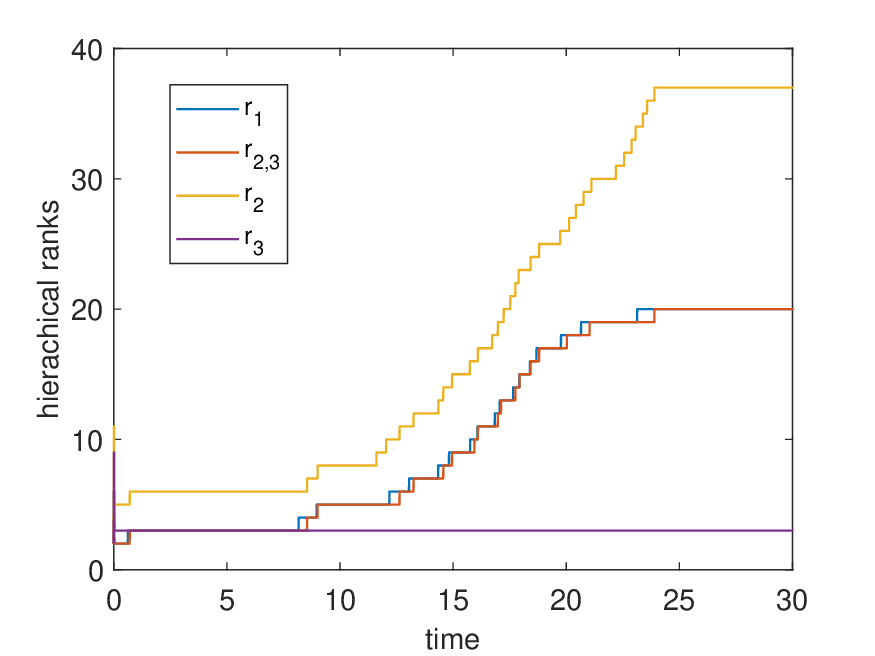}}
            \subfigure[]{\includegraphics[height=40mm]{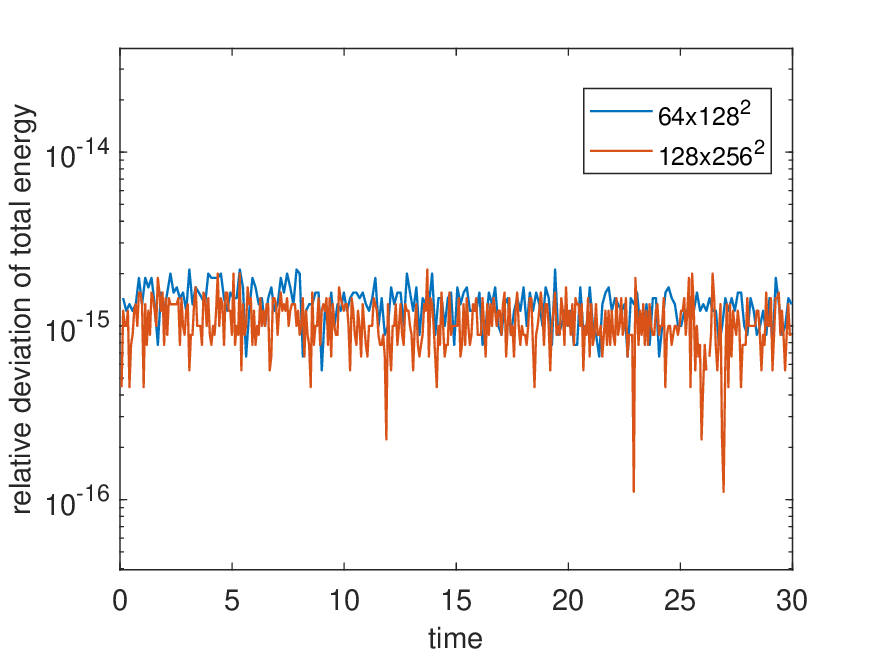}}
            \subfigure[]{\includegraphics[height=40mm]{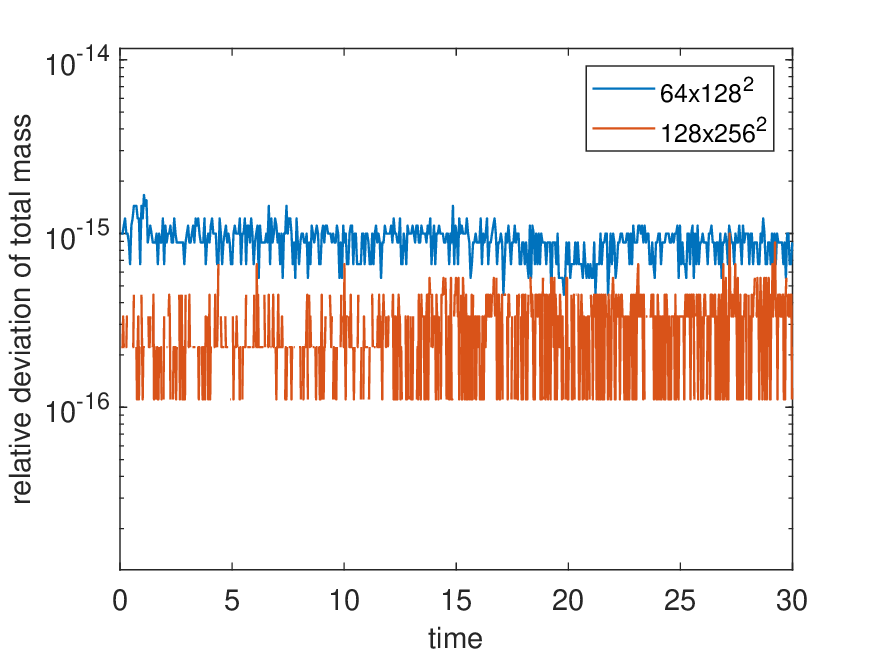}}
            
			\caption{Example \ref{ex:bump}. The time evolution of the electric energy for $E_1$, $E_2$ and magnetic energy (a), the hierarchical ranks of the numerical solution tensor (b, c), relative deviation of the total mass (d), and total energy  (e).  $T=30$, $\varepsilon=10^{-4}$. In (b) and (c) the curve of $r_{2,3}$ is close to that of $r_1$.}
			\label{fig:8}
		\end{figure}

\begin{figure}
			\centering
    \subfigure[$T=5$]{\includegraphics[height=60mm]{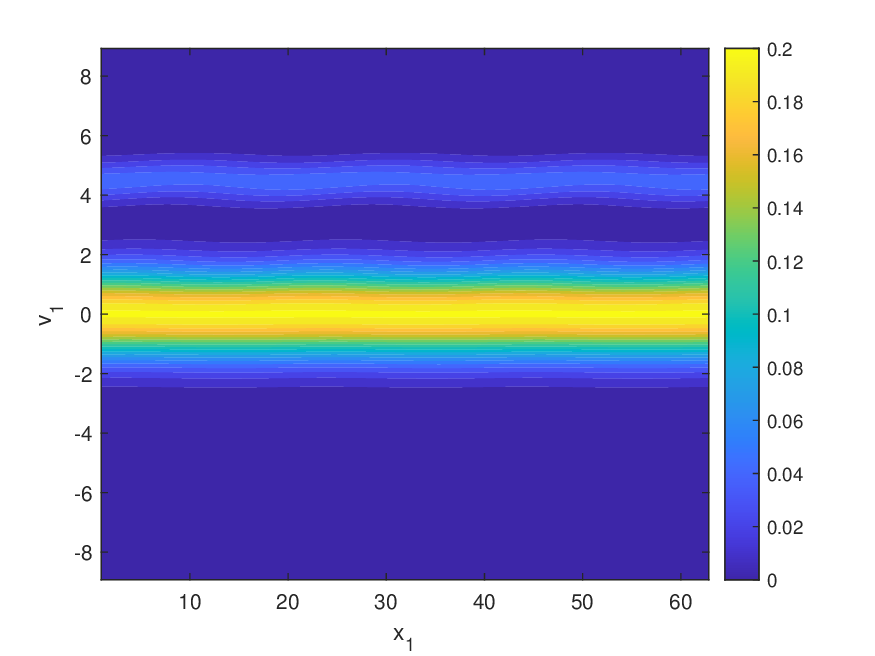}}
    \subfigure[$T=10$]{\includegraphics[height=60mm]{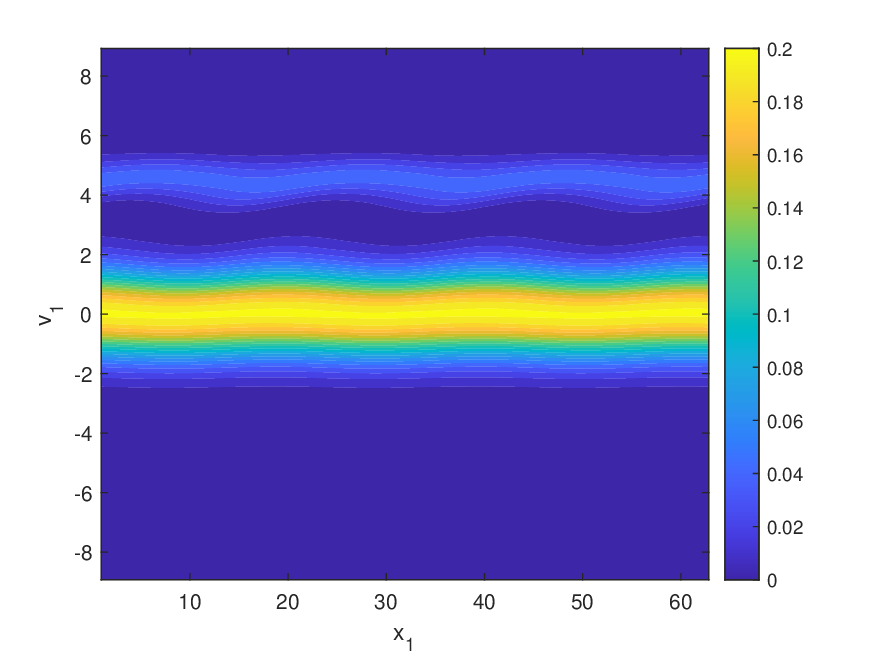}} 	
     \subfigure[$T=20$]{\includegraphics[height=60mm]{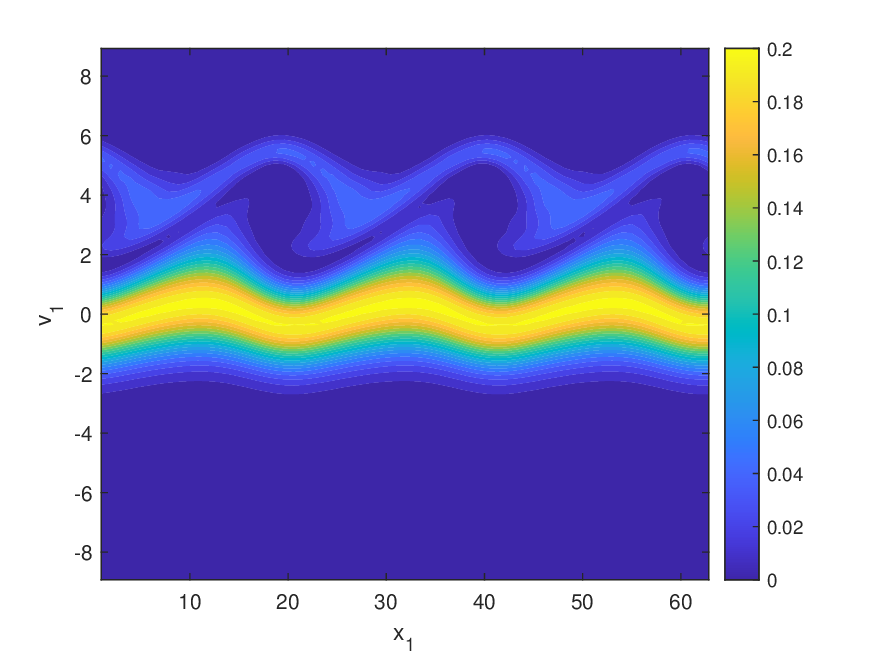}}
     \subfigure[$T=30$]{\includegraphics[height=60mm]{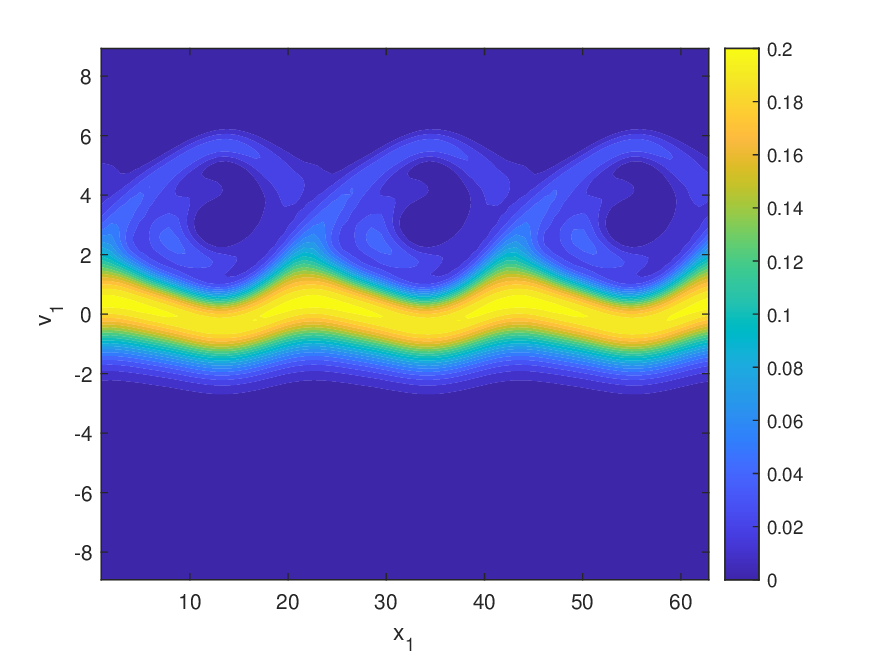}}
     \caption{Example \ref{ex:bump}. The contour plots of the 2D cuts of the solution  at $v_2=0.5$ for $N_x\times N_v^2 = 64 \times 128^2$ at $T=5,\,10,\,20,\,30$, and $\varepsilon = 10^{-4}$.}
     \label{fig:7}
\end{figure}

\end{exa}
\begin{exa}\label{ex:weibel}  In this example, we simulate the 1D2V Weibel instability \cite{cheng2014discontinuous,cheng2014energy}. The initial condition is 
    \[f(x_1,v_1,v_2,t=0) = \frac{1}{ \pi \beta}  e^{-v_{1}^2/\beta} \left( \delta e^{-(v_{2}-v_{0,1})^2 /\beta} + (1-\delta)e^{-(v_{2}+v_{0,2})^2 /\beta}\right), \]
   where $\beta = 0.01$.
The electric fields  $E_1$ and $E_2$ are set to zero initially, and the magnetic field at $t=0$ is chosen as
    \[B_3(x_1,t=0) = -\alpha \sin{(kx_1)},\] 
where $k = 0.2$ and $\alpha = 10^{-3}$. Similar to the two-stream instability, such perturbation in  $B_3$ will  start the Weibel instability.  The space domain here is $\Omega_{x_1}= [0, 2\pi/k]$ and the velocity domain is chosen as $\Omega_{v_1} \times \Omega_{v_2}=[-1.2, 1.2]^2$. We set $\varepsilon = 10^{-5}$ for truncation and two different sets of parameters will be considered,
\begin{align*}
    &\text{choice 1:} \hspace{3mm} \delta = 0.5,\, v_{0,1} = v_{0,2} = 0.3; \\
    &\text{choice 2:} \hspace{3mm} \delta = 1/6,\, v_{0,1} = 0.5,\, v_{0,2} = 0.1.
\end{align*}
For this example, we compare the non-conservative low rank method which does not employ any conservation correction techniques, including using
the conservative truncation algorithm and solving the macroscopic systems for the LoMaC property. In Figure \ref{fig:9}, we plot the contours of the 2D cuts for both methods at $x_1=0.05\pi$  for the two  parameter choices at at time $T=55,\,82$. The results produced by both low rank methods are consistent with  those reported in \cite{cheng2014discontinuous,cheng2014energy}, indicating that the proposed methods can correctly capture the dynamics of the Weibel instability. In  Figure \ref{fig:10}-\ref{fig:13}, we report the time histories of the electric and magnetic energy, numerical ranks and relative deviation of total mass and energy of the solutions computed by the non-conservative and conservative methods for comparison. It is observed that the evolution of the electric and magnetic energy agrees with those in \cite{cheng2014discontinuous}. Furthermore, the methods adapt the representation ranks  to efficiently resolve the  solution structures.
Lastly, the LoMaC low rank method is able to conserve the total mass and energy up to the machine precision, while the non-conservative method can conserve the invariants up to the truncation threshold.

\begin{figure}
			\centering
     \subfigure[Para. Choice 1: $T=55$]{\includegraphics[height=45mm]{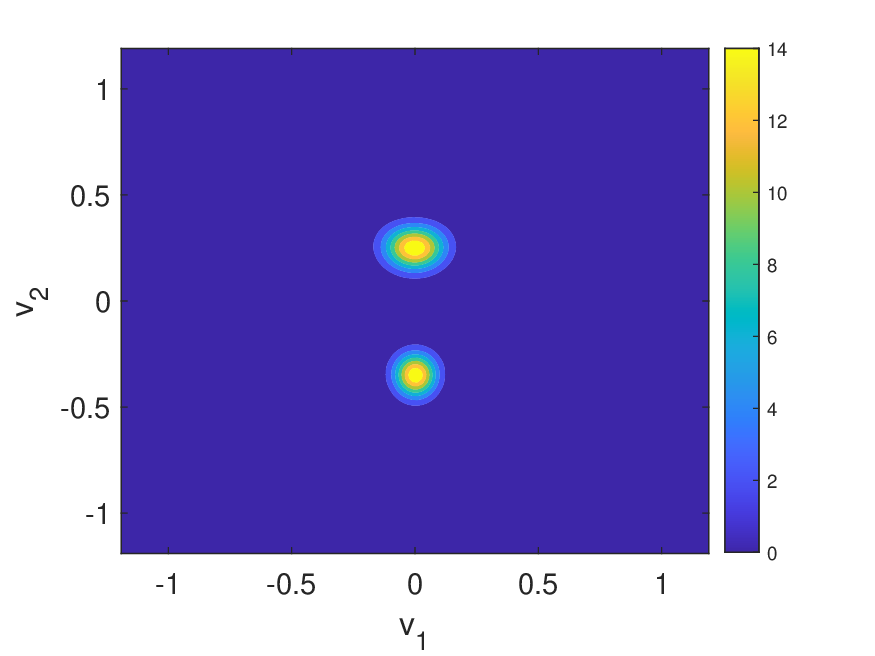}}
     \subfigure[Para. Choice 2: $T=55$]{\includegraphics[height=45mm]{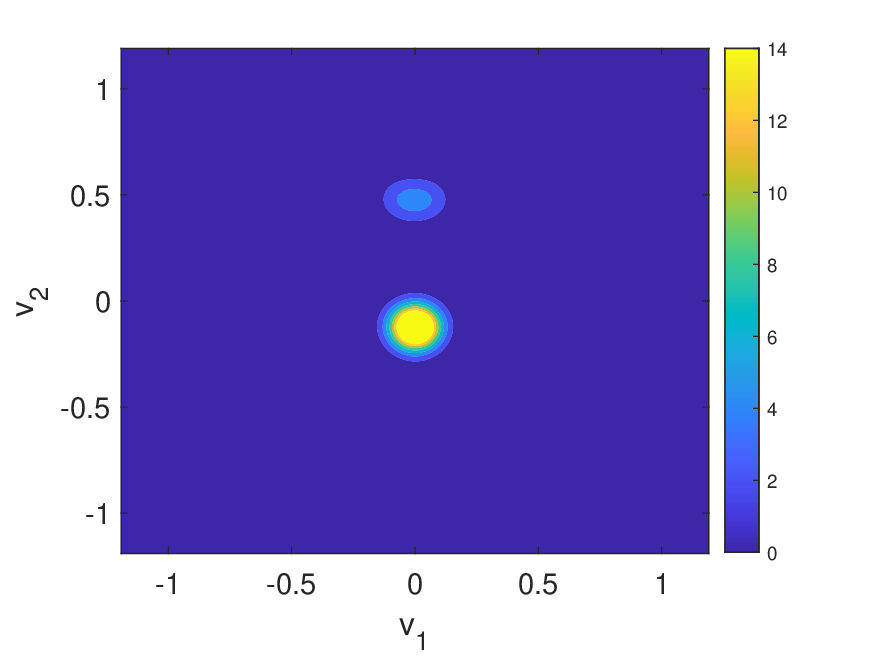}}
     \subfigure[Para. Choice 1: $T=55$]{\includegraphics[height=45mm]{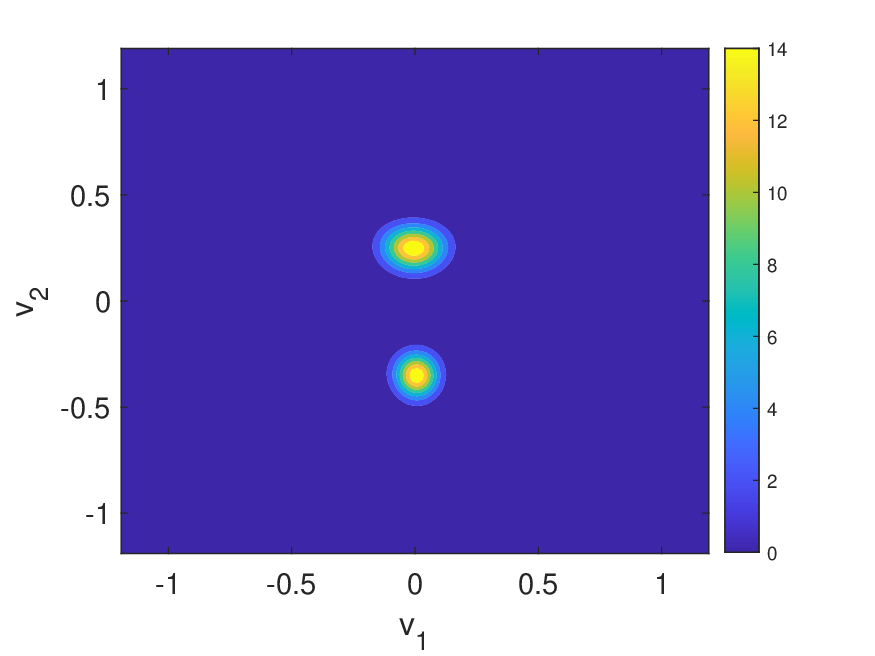}}
     \subfigure[Para. Choice 2: $T=55$]{\includegraphics[height=45mm]{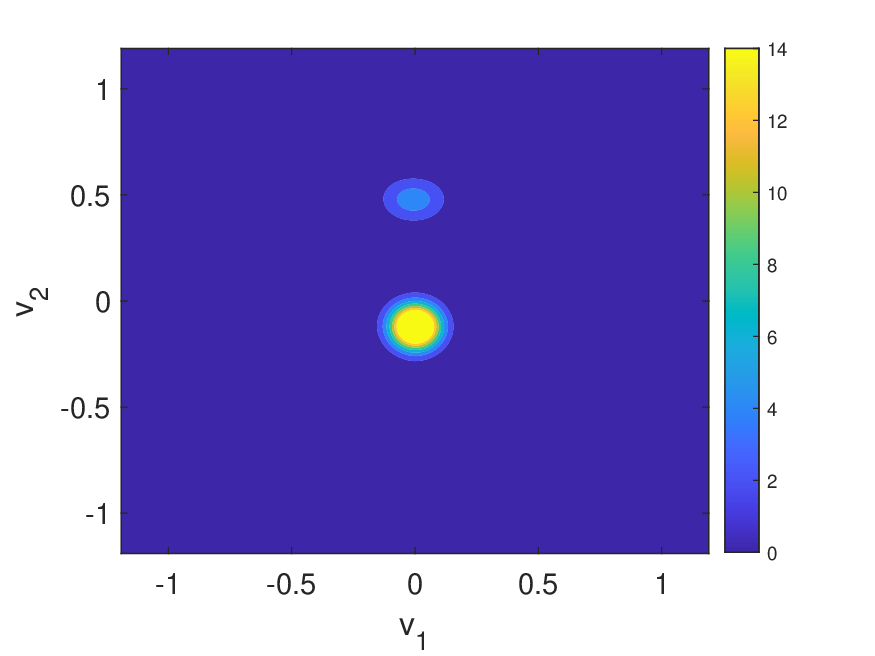}}
      \subfigure[Para. Choice 1: $T=82$]{\includegraphics[height=45mm]{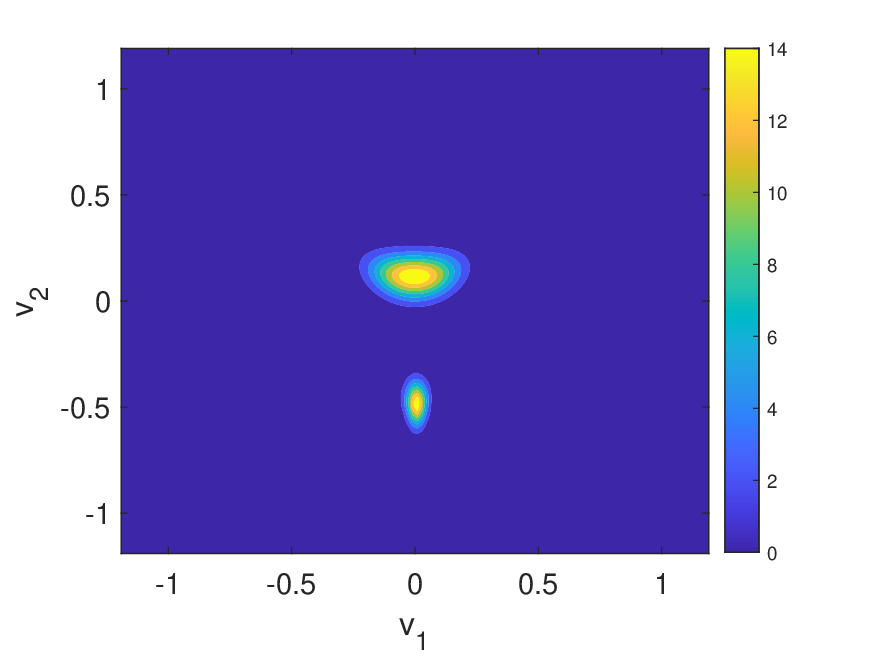}}
     \subfigure[Para. Choice 2: $T=82$]{\includegraphics[height=45mm]{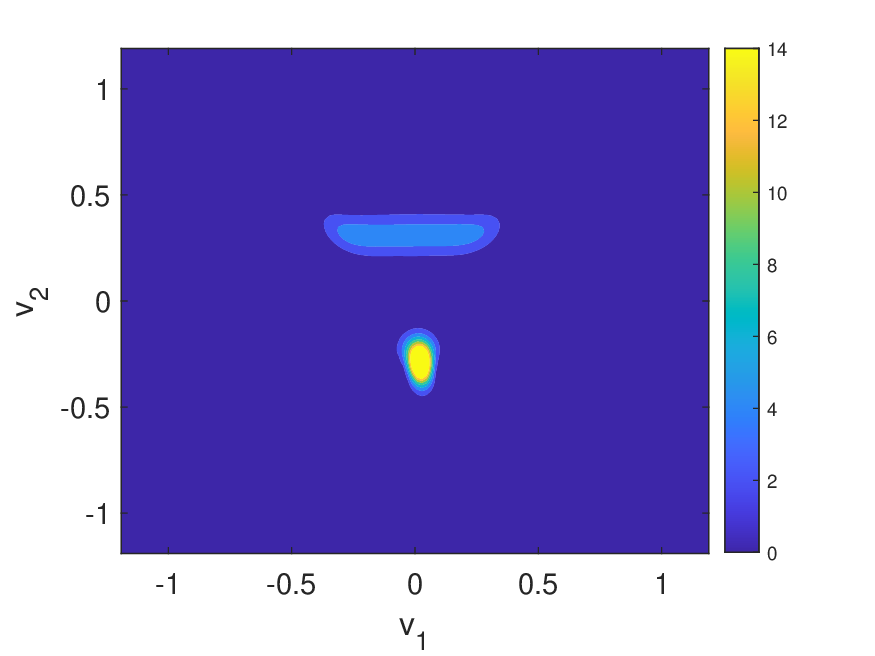}}
     \subfigure[Para. Choice 1: $T=82$]{\includegraphics[height=45mm]{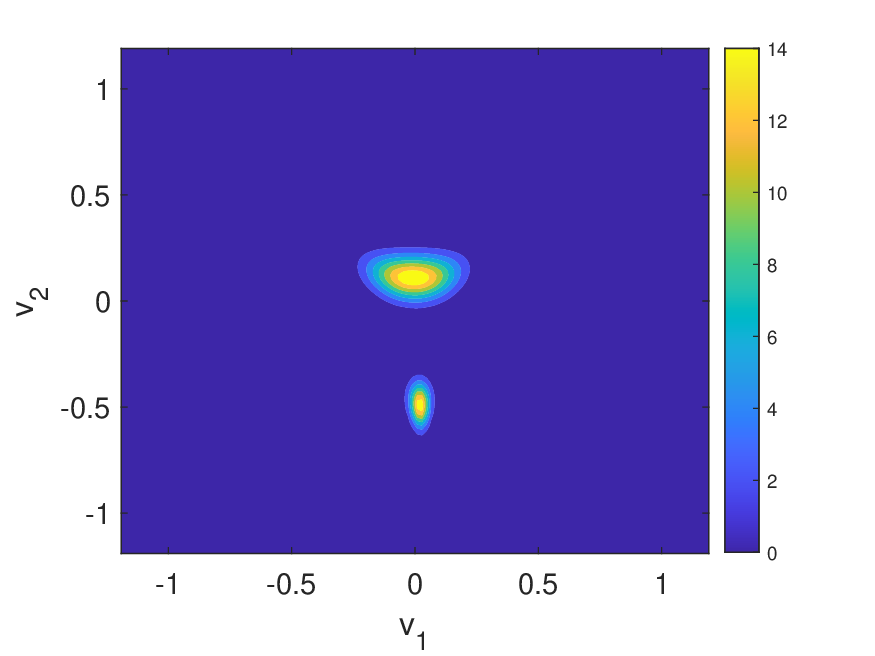}}
     \subfigure[Para. Choice 2: $T=82$]{\includegraphics[height=45mm]{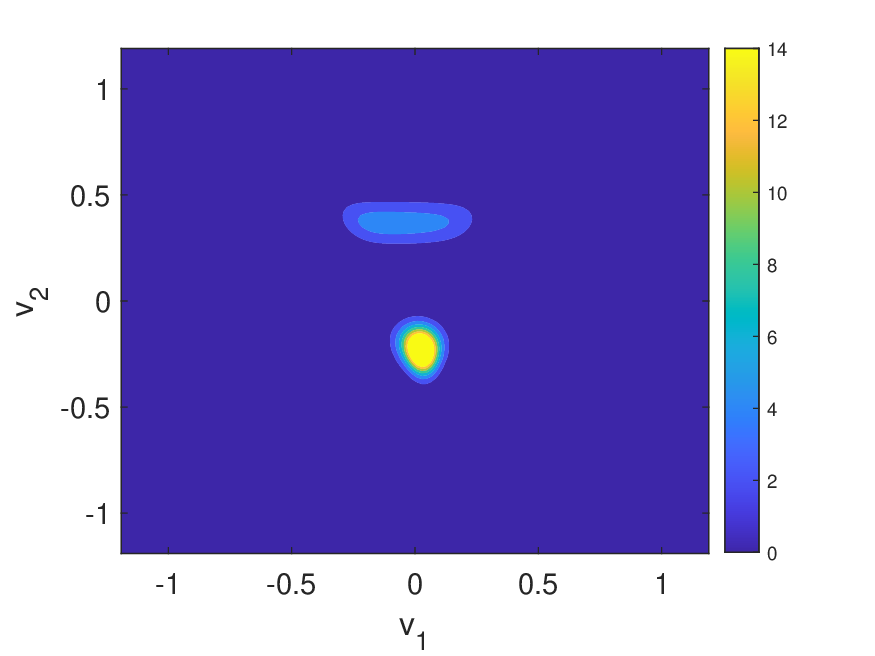}}
     \caption{ Example \ref{ex:weibel}. 2D contour plots at $x_1=0.05\pi$ for $N_x\times N_v^2 = 64 \times 128^2$ and $ \varepsilon = 10^{-5}$. Non-conservative method (a, b, e, f) and LoMaC method (c, d, g, h). }  
	\label{fig:9}
\end{figure}

  \begin{figure}
			\centering
     \subfigure[]{\includegraphics[height=40mm]{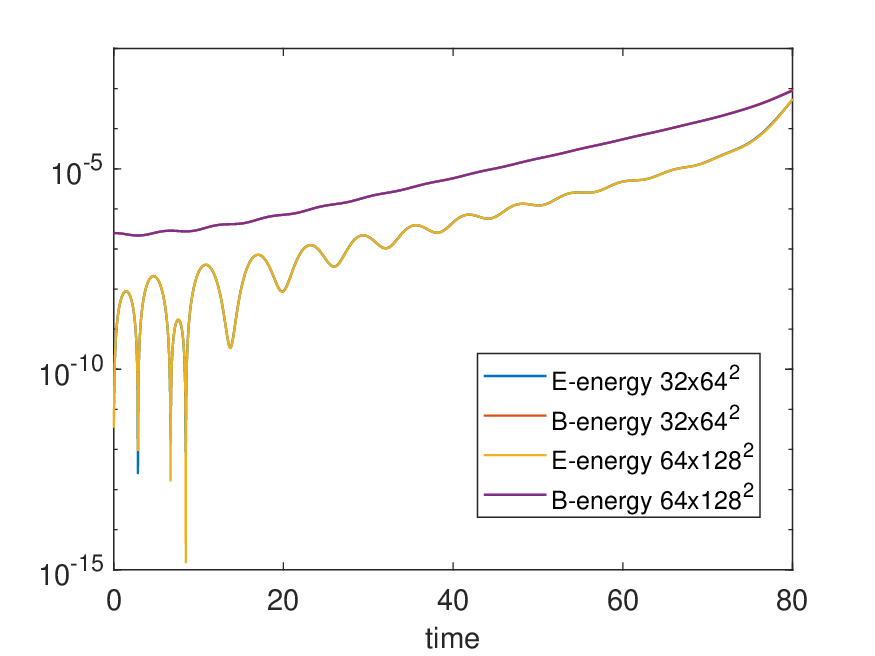}}
       \subfigure[$32 \times 64^2$]
     {\includegraphics[height=40mm]{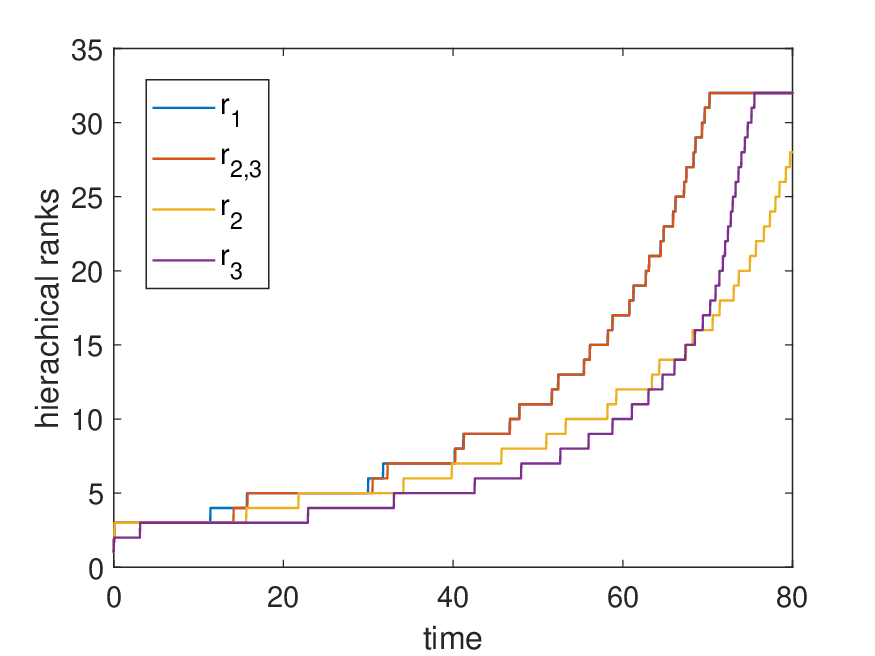}}
      \subfigure[$64 \times 128^2$]{\includegraphics[height=40mm]{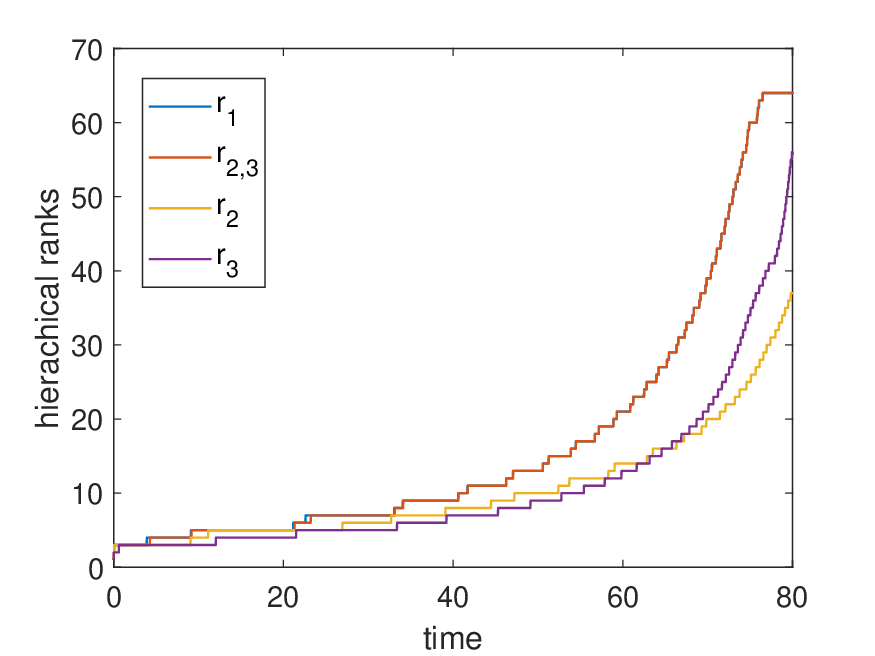}}
     \subfigure[]{\includegraphics[height=40mm]{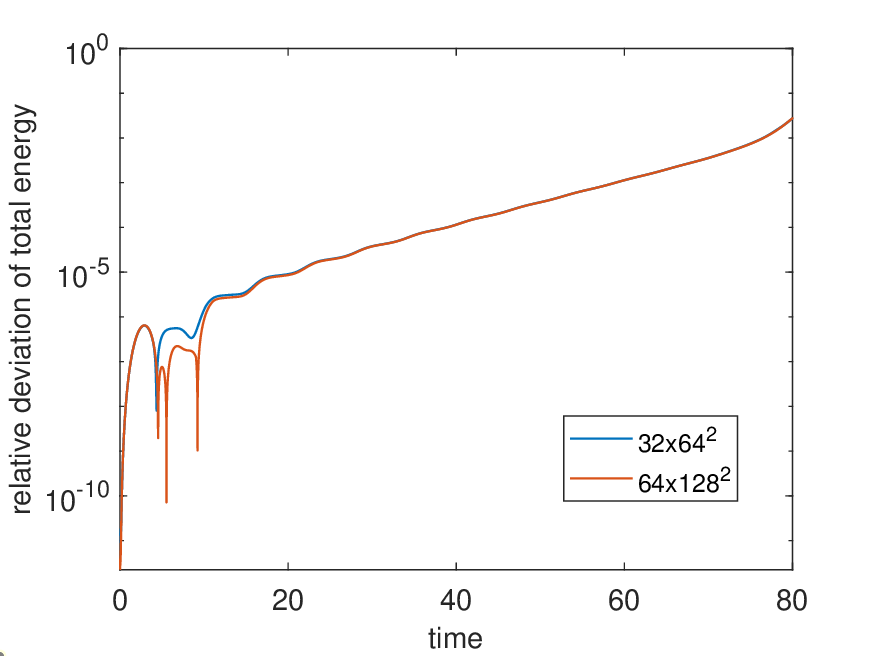}}
     \subfigure[]{\includegraphics[height=40mm]{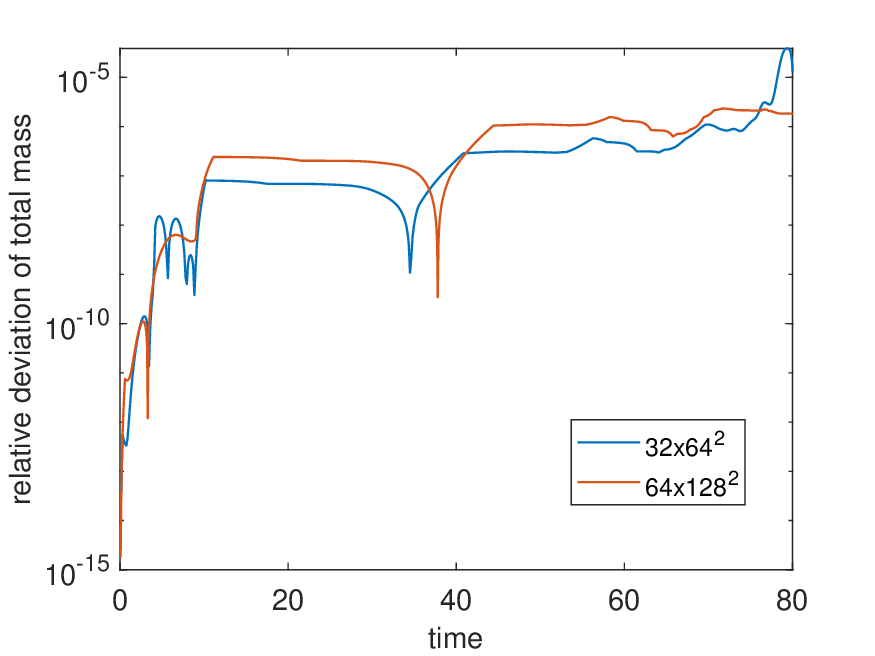}}
	\caption{Example  \ref{ex:weibel}. Parameter choice 1. The time evolution of the electric and magnetic energy for non-conservative method (a), the rank of the numerical solutions (b,c), relative
deviation of total energy (d), and total mass (e).  $ \varepsilon = 10^{-5}$.}
	\label{fig:10}
	\end{figure}	

\begin{figure}
			\centering
     \subfigure[]{\includegraphics[height=40mm]{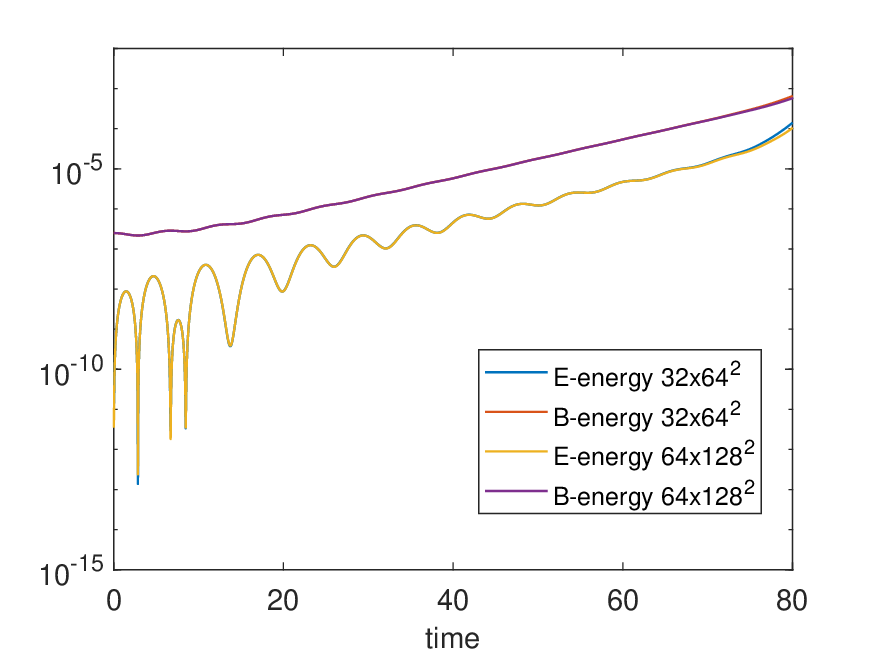}}
       \subfigure[$32 \times 64^2$]
     {\includegraphics[height=40mm]{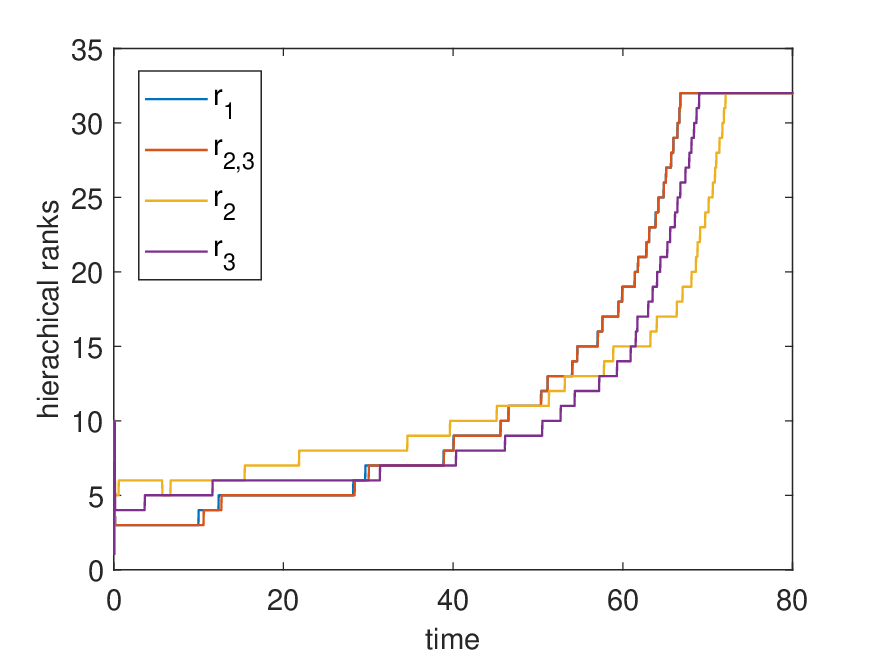}}
      \subfigure[$64 \times 128^2$]{\includegraphics[height=40mm]{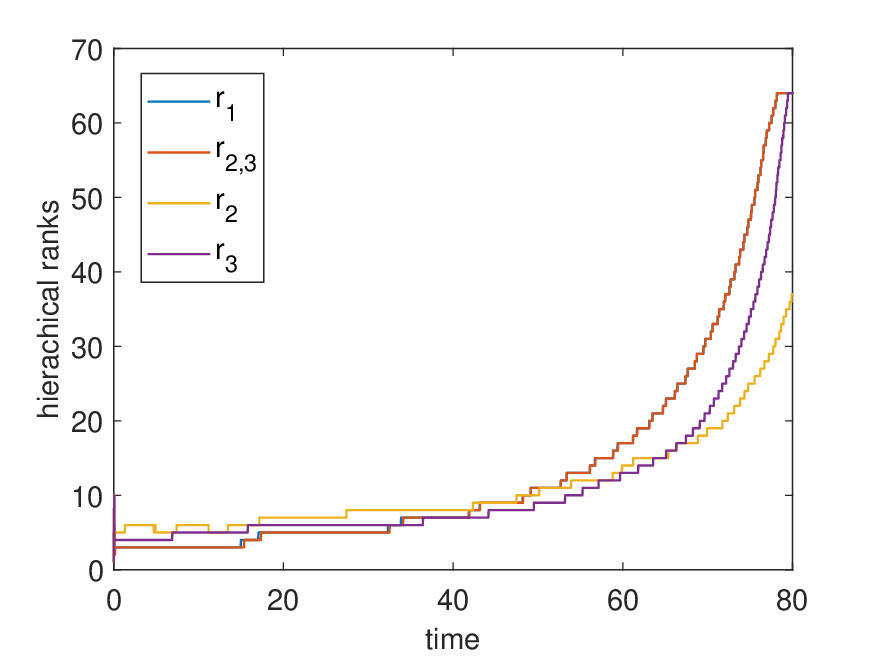}}
     \subfigure[]{\includegraphics[height=40mm]{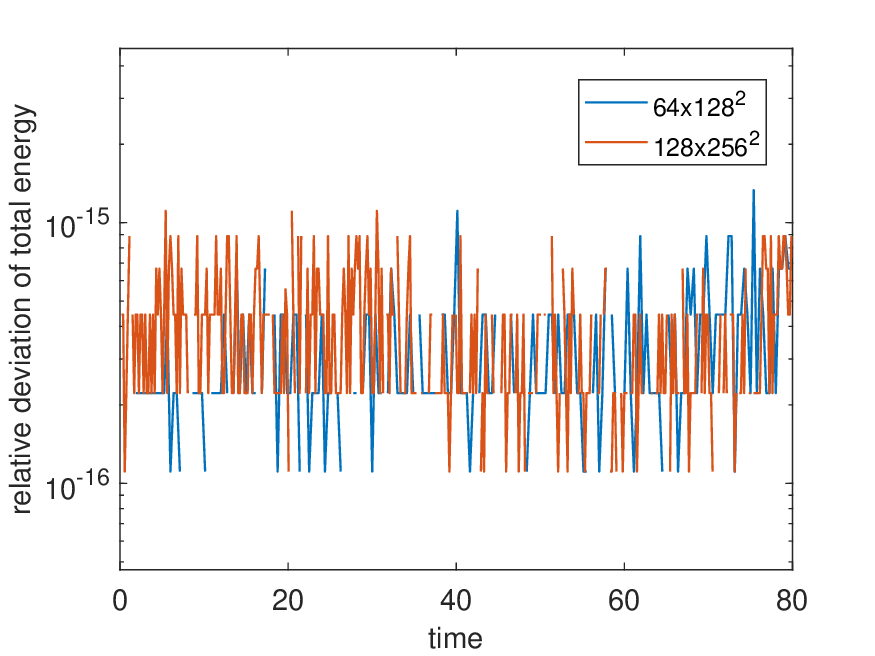}}
     \subfigure[]{\includegraphics[height=40mm]{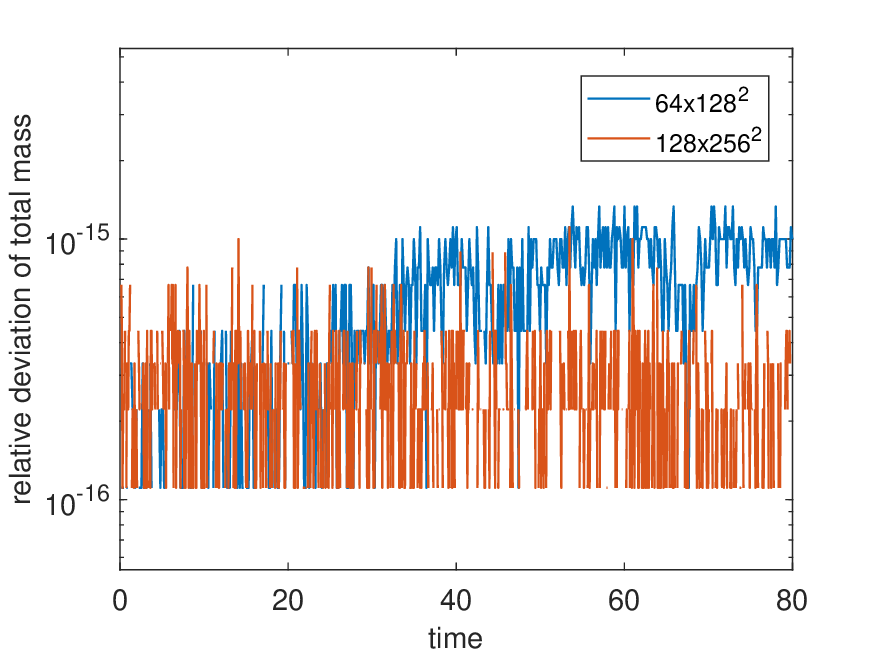}}
	\caption{Example  \ref{ex:weibel}. Parameter choice 1. The time evolution of the electric and magnetic energy for conservative method (a), the rank of the numerical solutions (b,c), relative
deviation of total energy (d), and total mass (e). $ \varepsilon = 10^{-5}$.}
	\label{fig:11}
	\end{figure}	

\begin{figure}
			\centering
     \subfigure[]{\includegraphics[height=40mm]{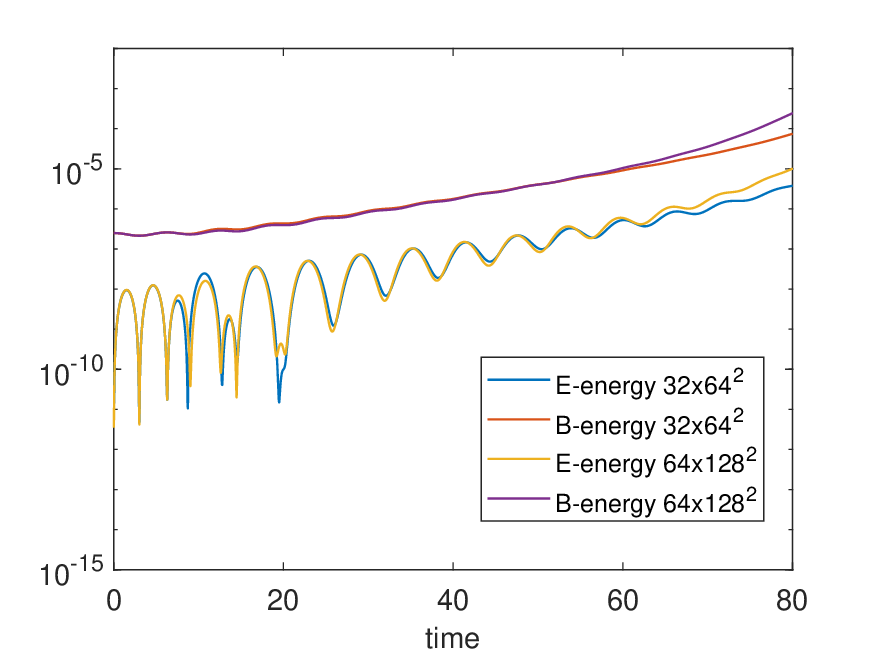}}
       \subfigure[$32 \times 64^2$]
     {\includegraphics[height=40mm]{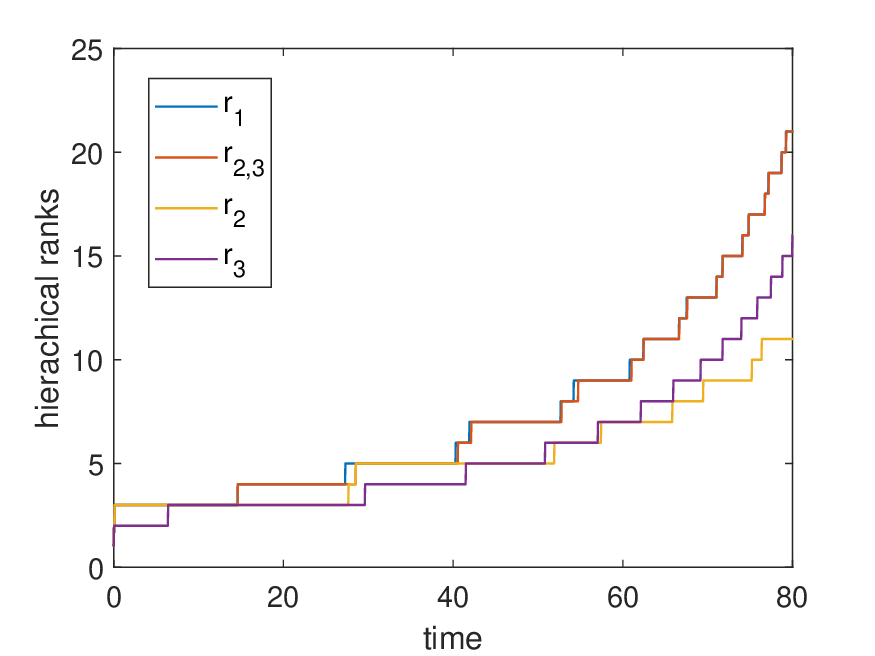}}
      \subfigure[$64 \times 128^2$]{\includegraphics[height=40mm]{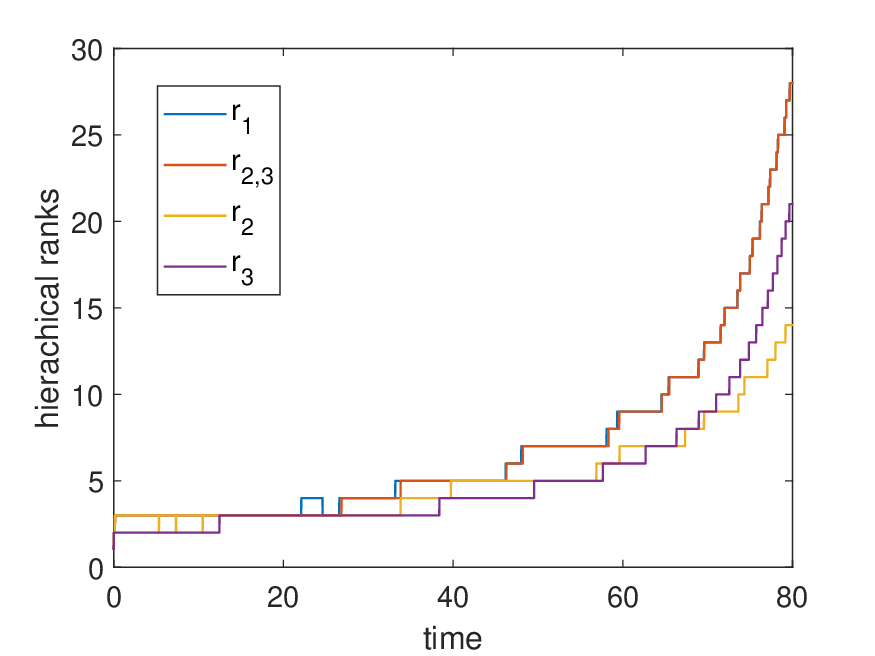}}
     \subfigure[]{\includegraphics[height=40mm]{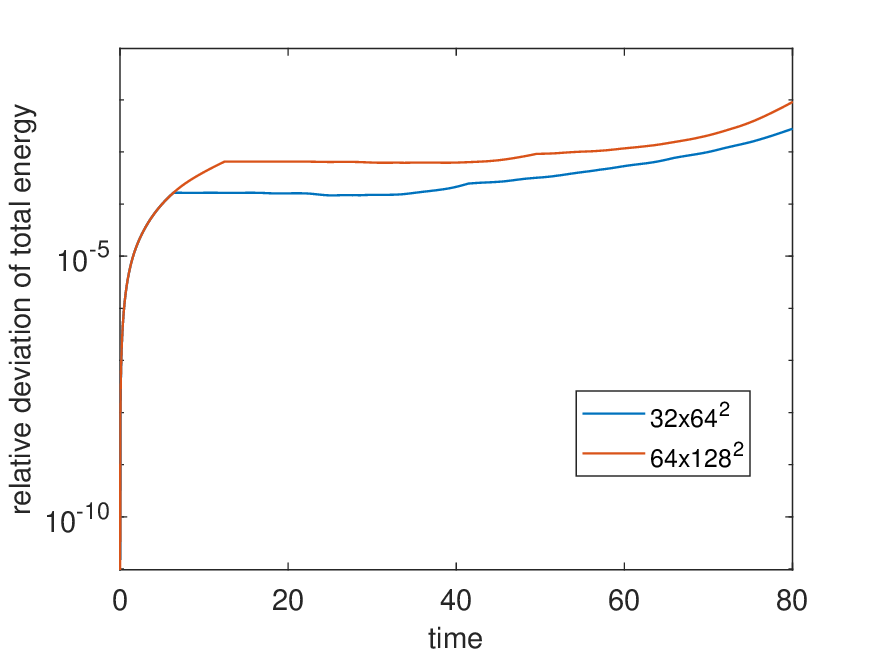}}
     \subfigure[]{\includegraphics[height=40mm]{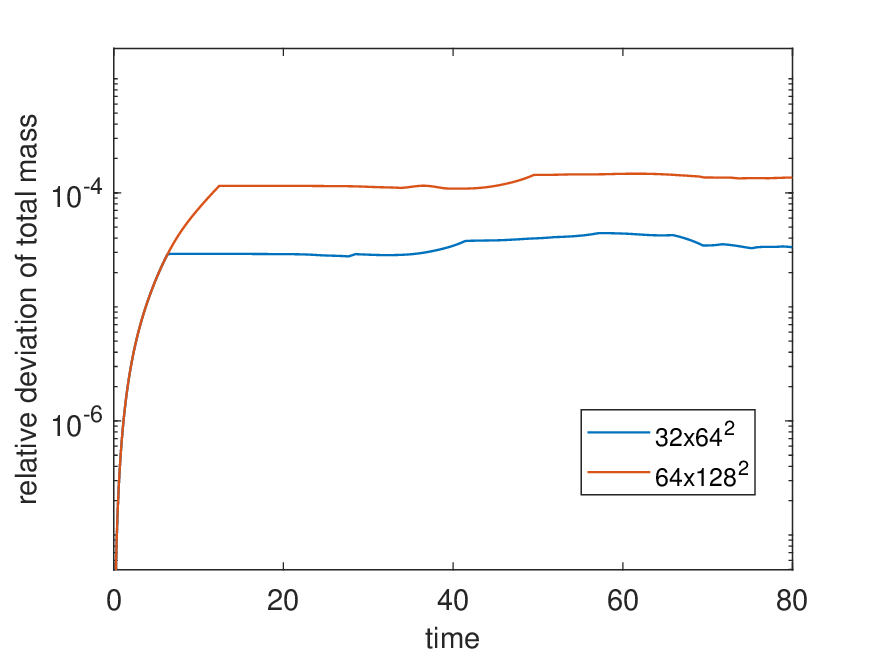}}
	\caption{Example  \ref{ex:weibel}. Parameter choice 2. The time evolution of the electric and magnetic energy for non-conservative method (a), the rank of the numerical solutions (b,c), relative
deviation of total energy (d), and total mass (e). $ \varepsilon = 10^{-5}$.}
	\label{fig:12}
	\end{figure}	

 \begin{figure}
			\centering
     \subfigure[]{\includegraphics[height=40mm]{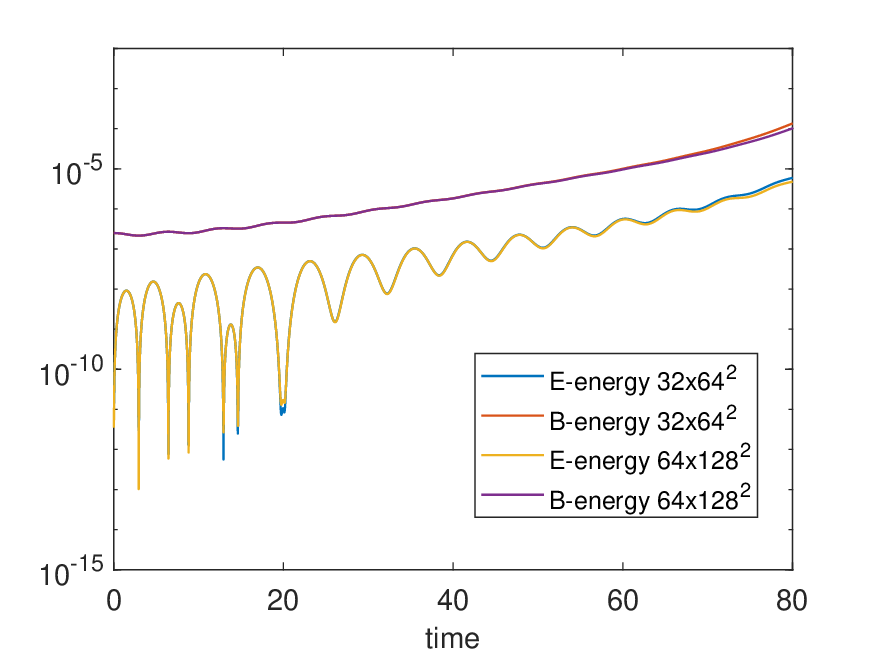}}
       \subfigure[$32 \times 64^2$]
     {\includegraphics[height=40mm]{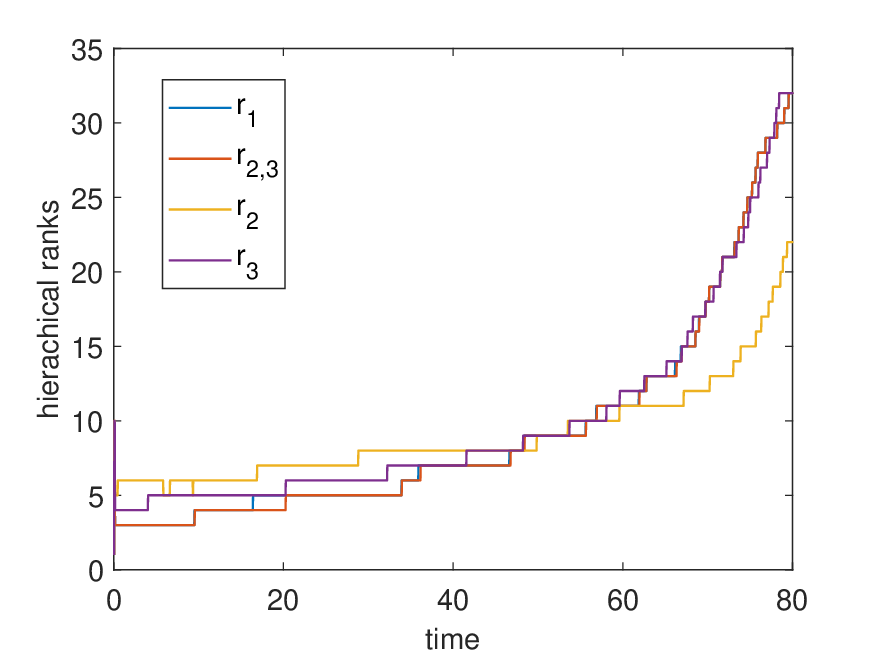}}
      \subfigure[$64 \times 128^2$]{\includegraphics[height=40mm]{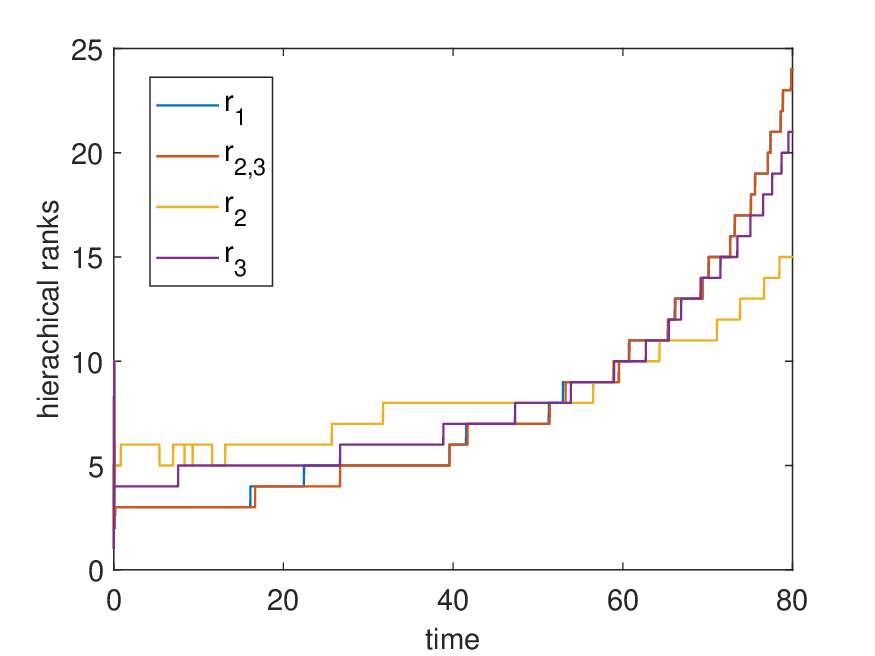}}
     \subfigure[]{\includegraphics[height=40mm]{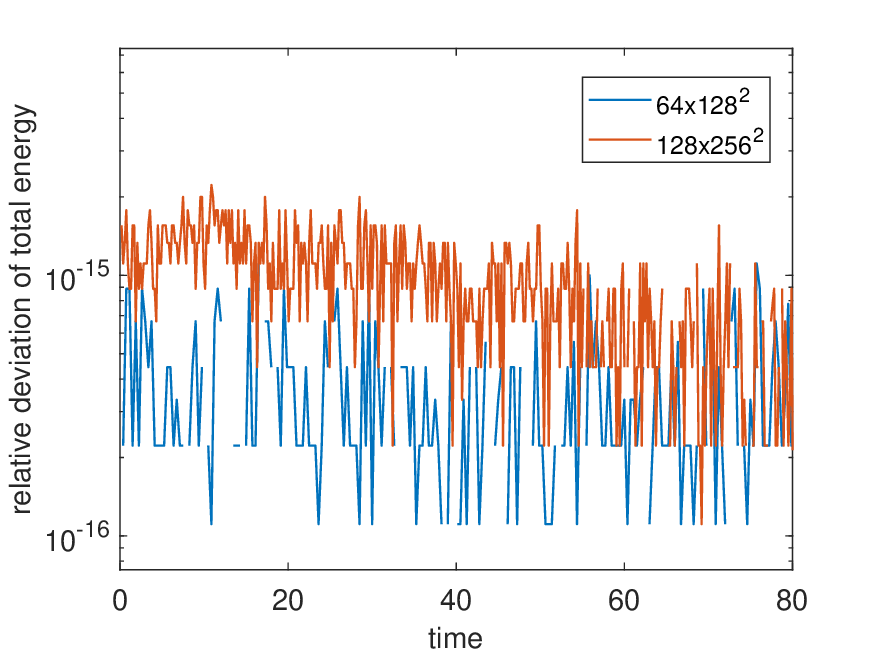}}
     \subfigure[]{\includegraphics[height=40mm]{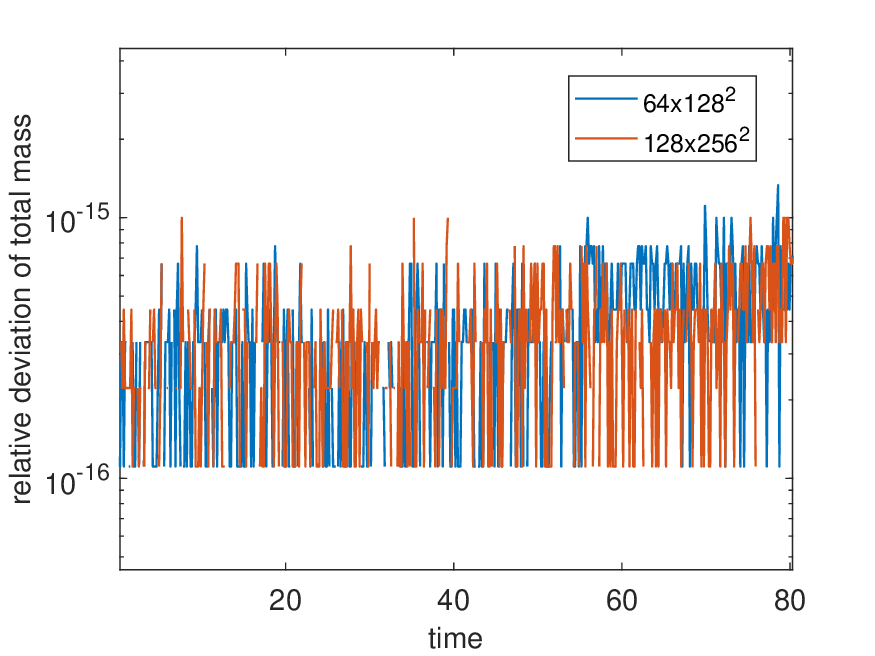}}
	\caption{Example  \ref{ex:weibel}. Parameter choice 2. The time evolution of the electric and magnetic energy for conservative method (a), the rank of the numerical solutions (b,c), relative
deviation of total energy (d), and total mass (e).  $ \varepsilon = 10^{-5}$.}
	\label{fig:13}
	\end{figure}

\end{exa}	

\begin{exa}
\label{ex:2d2v}
For the last example, we simulate the following the 2D2V version of the Landau-type problem with the initial condition 
    \[f(x_1,x_2,v_1,v_2,t=0) = \frac{1}{2 \pi} e^{-\frac{1}{2} (v_{1}^2+v_{2}^2)} \left(1+\alpha \left(\cos\left(kx_1\right) + \cos\left(kx_2\right)\right)\right),\]
where  $k = 0.4$ and $\alpha = 0.01$. The electric field $(E_1, E_2)$ is initialized by solving 2D Gauss’ law \eqref{eq:gauss2d}, and the magnetic field at $t=0$ is chosen as
    \[ B_3(x_1,x_2,t=0) = - E_1,\]
The spatial domain is $\Omega_{x_1}\times\Omega_{x_2}= [0, 2\pi/k]^2$, and the velocity domain is chosen as $\Omega_{v_1} \times \Omega_{v_2}=[-5, 5]^2$, similar to the 1D2V case. Note that we do not exploit the decomposition in $x_1$ and $x_2$ which will greatly facilitate the algorithm design to achieve LoMaC property. The truncation threshold is set to be  $\varepsilon = 10^{-5}$. In Figure \ref{fig:2d2v}, we report the time histories of the electric and magnetic energy,  numerical ranks and the time histories of the relative deviation of total mass and energy of the proposed LoMaC low rank solutions for two sets of meshes $N_x^2\times N_v^2 = 32^2 \times 64^2$ and $N_x^2\times N_v^2 = 64^2 \times 128^2$.  It is observed that the dynamics of electromagnetic energy is similar to the 1D2V case investigated in Example 3.1. That is the electric energy decays over time and displays oscillatory behavior. Again, the proposed LoMaC method is capable of conserving the total mass and energy up to the machine precision regardless of the mesh size. Furthermore, by adapting the hierarchical ranks of the solution tensor, the method can effectively capture the underlying dynamics.
 \begin{figure}
			\centering
     \subfigure[]{\includegraphics[height=40mm]{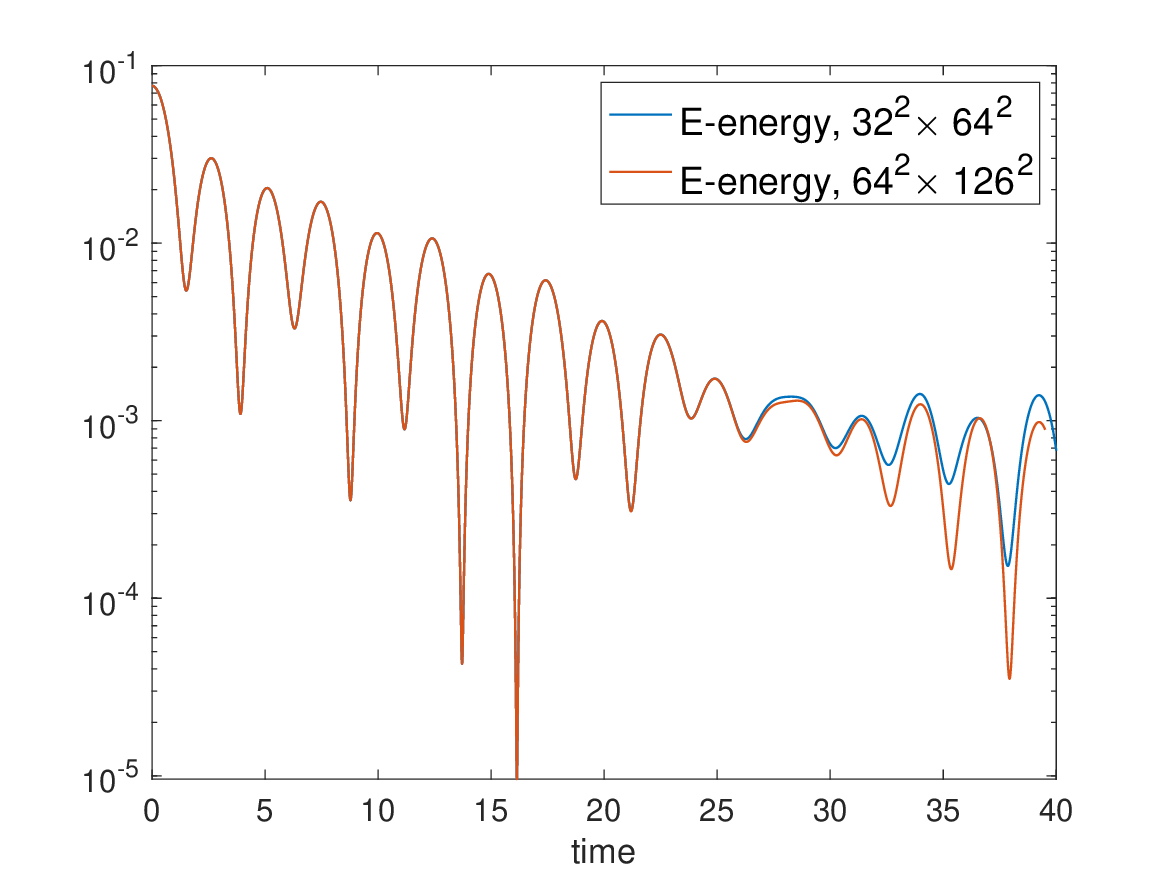}}
       \subfigure[]
     {\includegraphics[height=40mm]{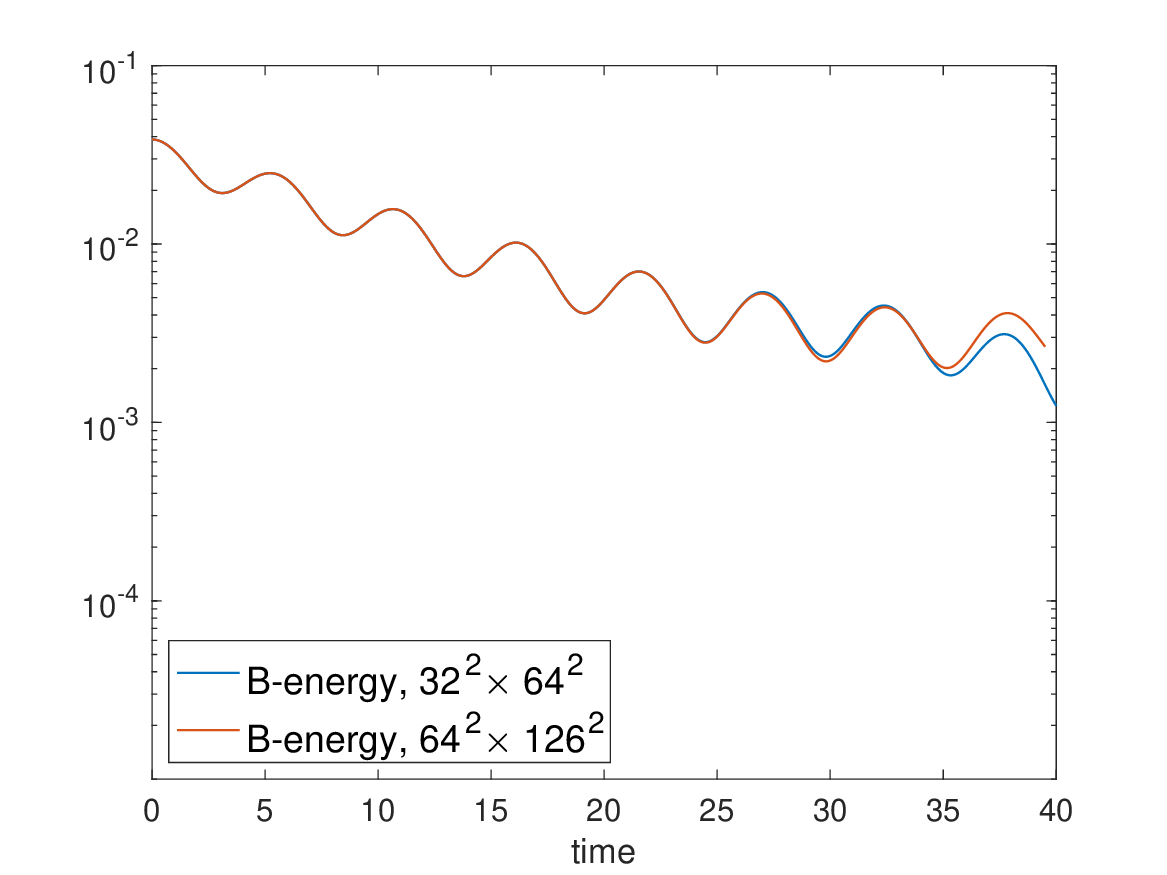}}
      \subfigure[]{\includegraphics[height=40mm]{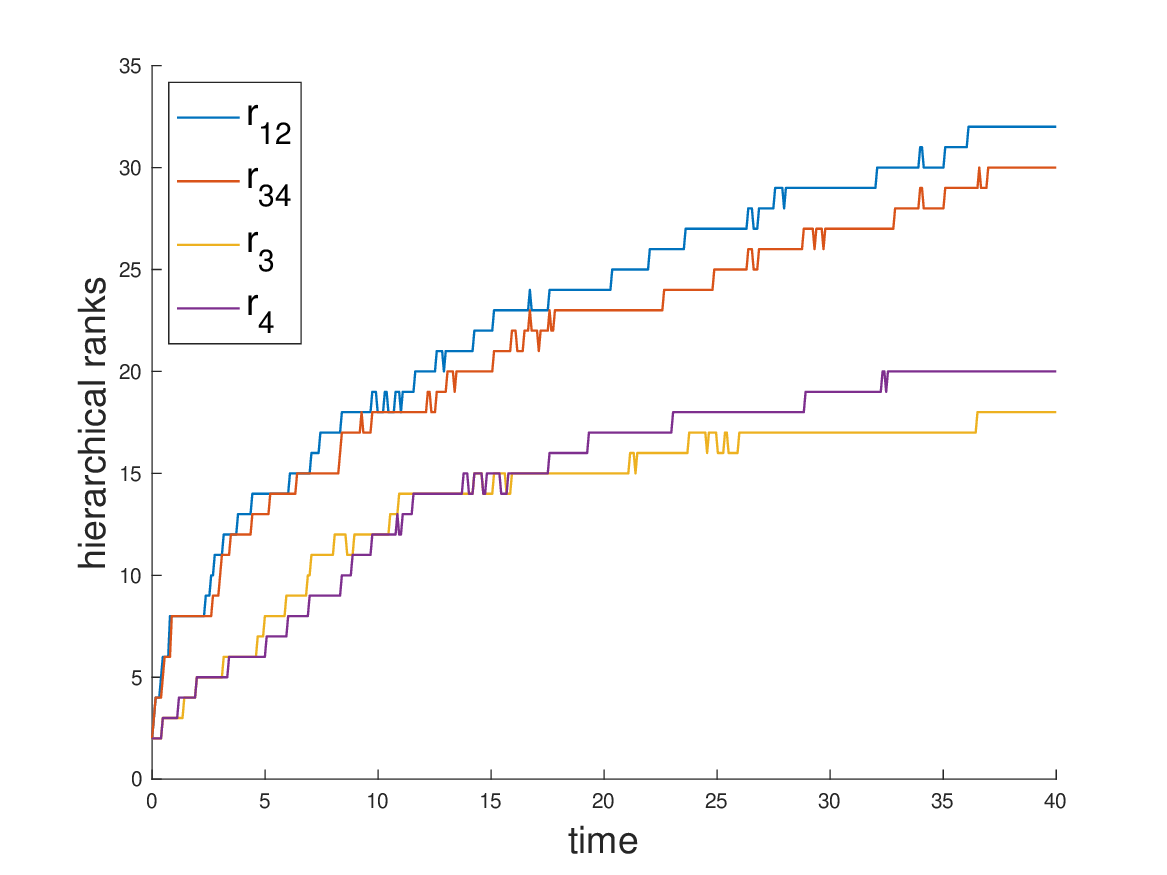}}
     \subfigure[]{\includegraphics[height=40mm]{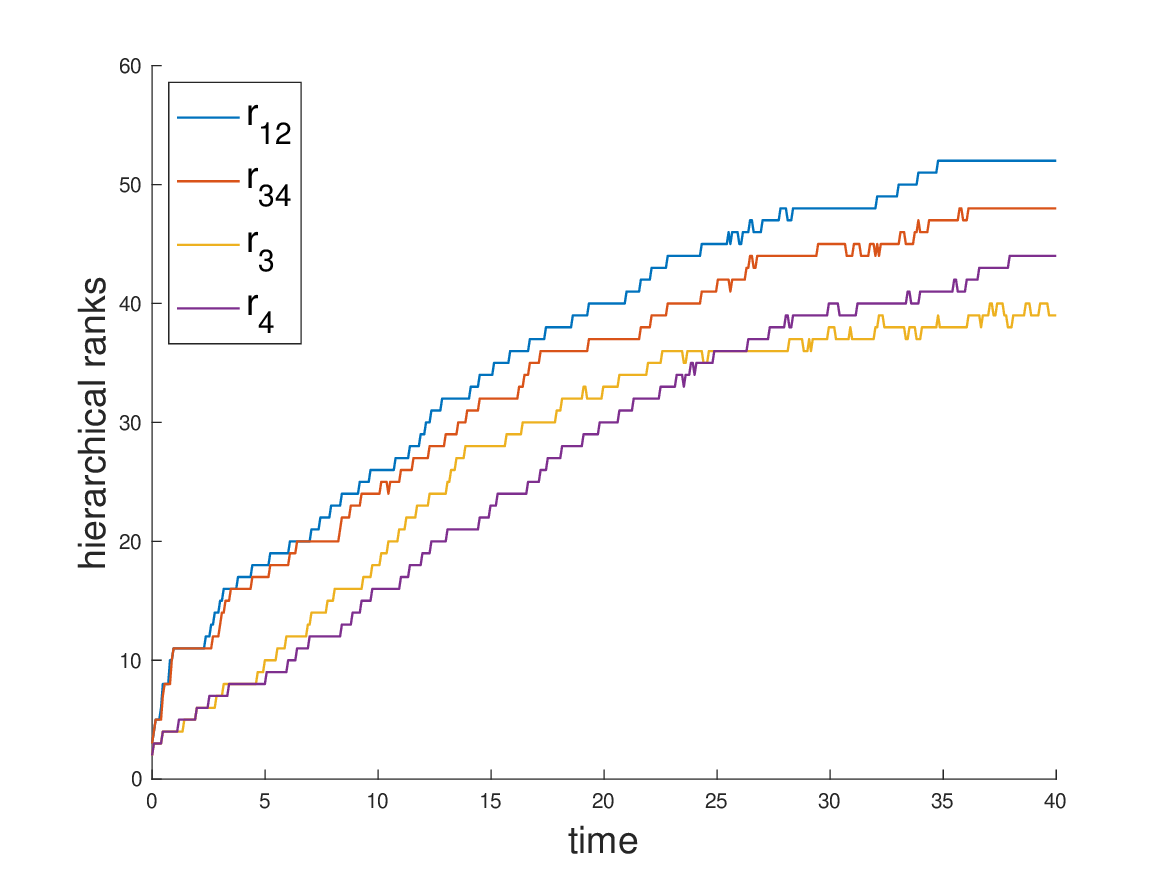}}
     \subfigure[]{\includegraphics[height=40mm]{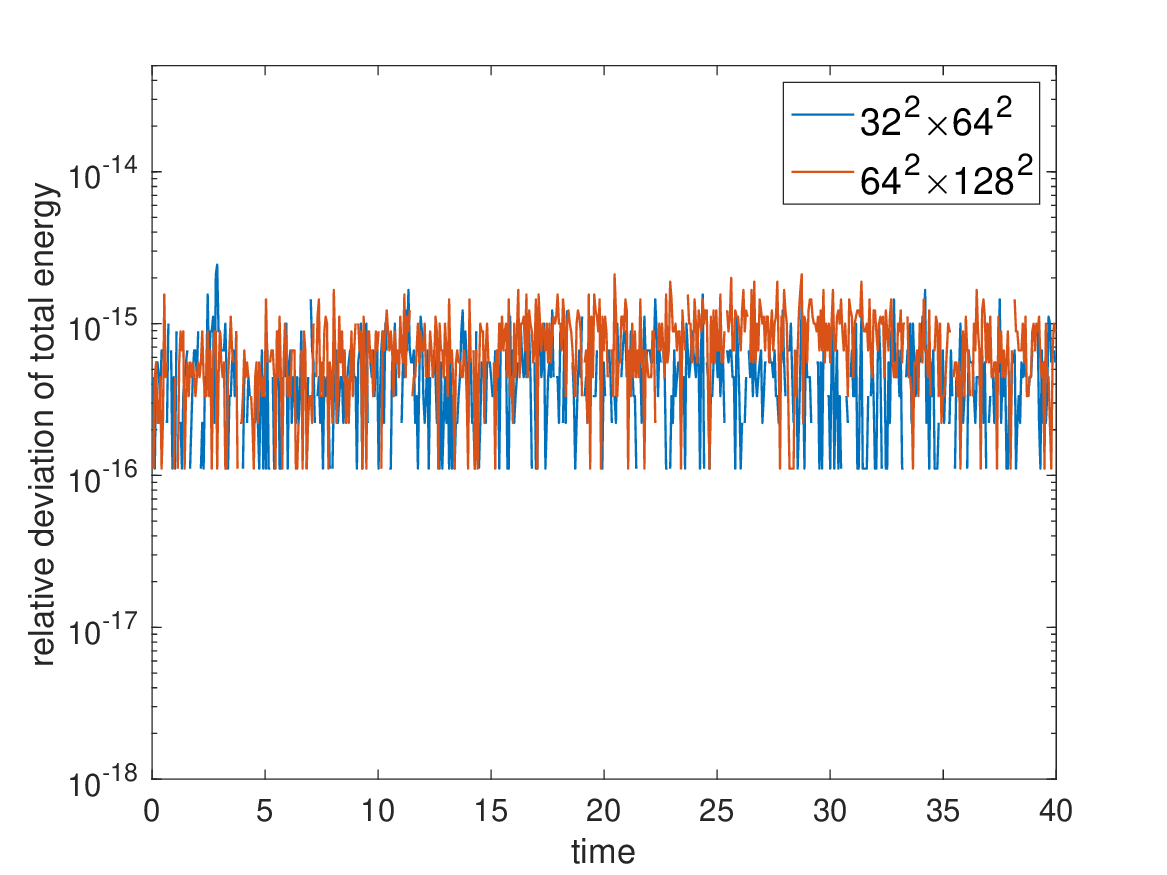}}
     \subfigure[]{\includegraphics[height=40mm]{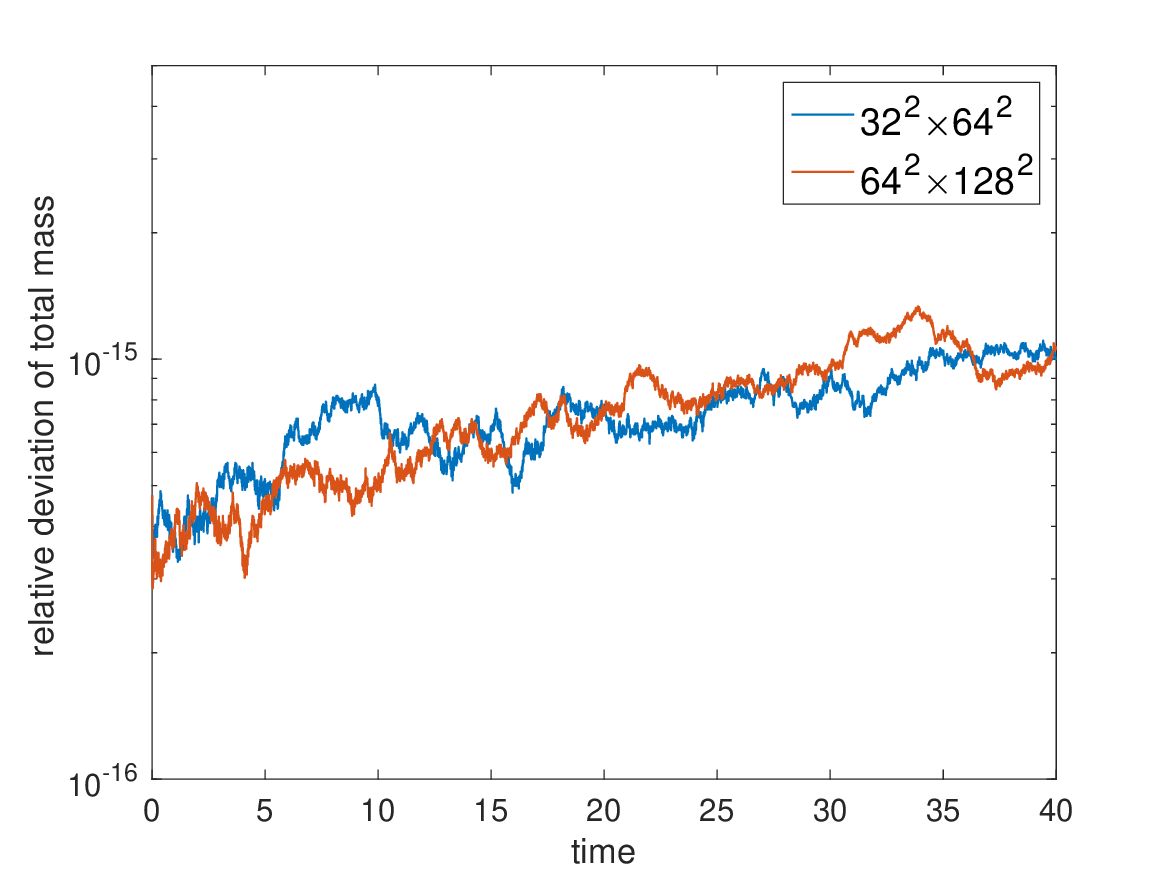}}
	\caption{Example  \ref{ex:2d2v}. The time evolution of the electric and magnetic energy for conservative method (a, b), the rank of the numerical solutions: $N_x^2\times N_v^2 = 32^2 \times 64^2$ (c), $N_x^2\times N_v^2 = 64^2 \times 128^2$ (d), relative
deviation of total energy (e), and total mass (f).  $ \varepsilon = 10^{-5}$.}
	\label{fig:2d2v}
	\end{figure}

\end{exa}

%% file: conclusion.tex
\section{Conclusion}
In this paper, we developed a low rank tensor method for simulating the Vlasov-Maxwell (VM) system. The method makes use of the hierarchical Tucker tensor format to express the kinetic solution in high dimensions. Together with the conservative truncation algorithm and working with  the macroscopic conservation laws,  the proposed method can achieve the  Local Macroscopic Conservative (LoMaC) property with explicit time marching. The future work includes the extension to collisional kinetic models and the development of positivity preserving technique under the proposed LoMaC low rank framework.